\documentclass[]{elsarticle}
\usepackage{indentfirst}
\usepackage{multirow}
\usepackage{tikz}
\usetikzlibrary{arrows,shapes}
\usepackage{fancybox}
\usetikzlibrary{shapes.multipart}
\usepackage{graphicx}
\usepackage{subcaption}
\usepackage{float}
\usepackage{pgfplots}
\pgfplotsset{width=10cm,compat=1.9} %1.14

\RequirePackage{amsmath}
\RequirePackage{amssymb}
\RequirePackage{amsthm}
\RequirePackage{relsize}
\RequirePackage{graphicx}

\newcommand{\Restrict}[2]{{\left.\kern-\nulldelimiterspace#1\vphantom{\big|}\right|_{#2}}}

\newcommand{\Stress}{\sigma}    % stress tensor
\newcommand{\bm}{\boldsymbol}

\newcommand{\bff}{{\bm f}}

\newcommand{\bfn}{{\bm n}}

\newcommand{\bft}{{\bm t}}
\newcommand{\bfu}{{\bm u}}
\newcommand{\bfv}{{\bm v}}
\newcommand{\bfw}{{\bm w}}
\newcommand{\bfx}{{\bm x}}

\newcommand{\uu}{{\bfu}}
\newcommand{\vv}{{\bfv}}
\newcommand{\ww}{{\bfw}}
\newcommand{\xx}{{\bfx}}

\newcommand{\bfsigma}{{\boldsymbol{\Stress}}}

\journal{CMAME }
\begin{document}
%---------------------------------------------------
% PREAMBLE PART (TITLE, AUTHORS, KEYWORDS, ABSTRACT)
\begin{frontmatter}
\title{
Stabilized finite element method for incompressible solid dynamics using an updated Lagrangian formulation
}

\cortext[cor1]{~Corresponding author}
\address[]{MINES~ParisTech, PSL~-~Research~University, CEMEF~-~Centre~for~material~forming,\\ CNRS~UMR~7635, CS~10207~rue~Claude~Daunesse, 06904~Sophia-Antipolis~Cedex, France.}

\author[]{R.~Nemer\corref{cor1}}
\ead{ramy.nemer@mines-paristech.fr}
\author[]{A.~Larcher}
\ead{aurelien.larcher@mines-paristech.fr}
\author[]{T.~Coupez}
\ead{thierry.coupez@mines-paristech.fr}
\author[]{E.~Hachem}
\ead{elie.hachem@mines-paristech.fr}

\begin{abstract}
This paper proposes a novel way to solve transient linear, and non-linear solid dynamics for compressible, nearly incompressible, and incompressible material in the updated Lagrangian framework for tetrahedral unstructured finite elements.
It consists of a mixed formulation in both displacement and pressure, where the momentum equation of the continuum  is complemented with a pressure equation that handles incompresibility inherently.
It is obtained through the deviatoric and volumetric split of the stress, that enables us to solve the problem in the incompressible limit.
The Varitaional Multi-Scale method (VMS) is developed based on the orthogonal decomposition of the variables, which damps out spurious pressure fields for piece wise linear tetrahedral elements.
Various numerical examples are presented to assess the robustness, accuracy and capabilities of our scheme in bending dominated problems, and for complex geometries.
\end{abstract}

\end{frontmatter}
%
% END OF PREAMBULE
%---------------------------------------------------
% \documentclass[a4paper]{elsarticle}

\pagenumbering{arabic} 
\setcounter{page}{1}

\section{Introduction}
\indent

The need of a solid solver that can handle complex geometry is at its highest peaks.
Whether it be in everyday life applications, such as the behavior of electrical wires \cite{TAGHIPOUR2018211}, or the complex electro-elasticity behavior of components \cite{POYA201875}, to the understanding the elastic shock in solids \cite{8851076}.
Unconventional shapes, which result from complex algorithms, such as shape optimization, need to be tested.
In addition, bio-medical, and bio-mechanical applications usually include complex geometries based on the human organs.
This can also be of high of interest in the field of AI, which is generating atypical geometries through coupling with aforementioned fields.
Furthermore, a high range of materials can be considered incompressible or nearly incompressible.
Ranging from certain polymers that do not undergo high volume changes, to biological tissues.
\\
\\
\indent
A finite element formulation in which the displacement field is the unknown, and all other physical quantities are obtained using post-processing methodology  is typically used in solid dynamics \cite{470a3479a2a5460e9ad672060ea67376} \cite{10.5555/3226237.3226350}.
This method performs poorly near the incompressibility limit.
Locking, spurious pressure fields, and poor performance in bending related applications are some of the shortcomings of the preceding formulation \cite{hughes2012finite}.
\\
\\
\indent
This subject has seen its fair share of developments, aiming to damp or eliminate the above stated limitations.
Selective and reduced integration's methods, such as the B-bar \cite{hughes2012finite} \cite{10.1115/1.3423994} \cite{doi:10.1002/nme.1620150914}, the F-bar \cite{doi:10.1002/nme.1187} \cite{doi:10.1002/cnm.1640111109} \cite{DESOUZANETO19963277} \cite{MASUD2013359}, or the mean dilatation finite element methods \cite{nagtegaal1974numerically} are used for their ease of implementation.
By reducing the order of incompressibility at quadrature points, these methods circumvent the numerical instabilities of the inf-suf or Ladyzhenskaya-Babuska-Brezzi(LBB) conditions \cite{babuvska1971error}. 
These strategies have proven to be accurate for structured quadrilateral, and hexahedral meshes.
It is noteworthy to mention that automated grid generation for hexahedral elements for complex geometries, is very costly, and requires more time than the actual computing time.
This is however surpassed when using unstructured tetrahedral elements.
\\
\\
\indent
When considering the case of static, incompressbile elasticity, we obtain an elliptic equation, similar to that of the Stokes problem in fluid--mechanics, while the transient case or elastodynamics leads to a hyperbolic equation.
Given the similarities in the equations, it is natural to extend the mixed/coupled velocity/pressure formulation of the stokes problem \cite{HUGHES198685}, to the mixed displacement/pressure problem of the static elastic case \cite{10.1007/BF01395881}.
This extension acts as a bridge for the different, already implemented methods in fluid mechanics, to the solid mechanics field.
\\
\\
\indent
In \cite{chiumenti2002stabilized}, an incompressible steady-state linear elastic material was modeled using the mixed formulation in displacement/pressure using the Orthogonal Sub Scale method \cite{hughes1998variational}.
This work showed the capabilities of the mixed formulation (displacement/pressure) using the OSS method in the incompressible limit.
Also in \cite{cervera2010mixed} \cite{cervera2010mixed2}, an incompressible non-linear material was also modeled using the Orthogonal Sub Scale method.
These works, shows the capabilities of a mixed formulation using strain/displacement or stress/displacement formulation. 
There is a compromise however between computational cost and accuracy.
In \cite{article}, a three field (displacement/pressure/strain) formulation was tested and showed to be effective and accurate in the near incompressible limit.
\\
\\
\indent
Most of the former formulations were developed for the steady-state solution, and as mentioned before, transient elastodynamics convert the parabolic problem to an elliptic problem.
This is due to the second order derivative of displacement of the momentum equation.
This problem was addressed in different works, trying to circumvent this issue.
\\
\\
\indent
Some of these work include \cite{Bonet1}\cite{LEE201313}\cite{2014JCoPh.259..672A}\cite{leefirst}, where a finite-strain non-linear solid dynamics model is based on a new first-order (mixed) form of the equations in the Lagrangian framework.
The proposed methodology consists of adding an additional variable, which is the deformation tensor $\textbf{F}$, and Lagrange multipliers for the conservation of angular momentum if needed.
The results obtained are second order accurate in stress.
Moreover, in the incompressible limit, and bending dominated problems, an additional variable was introduced, which is the Jacobian determinant of the deformation gradient $J$ \cite{gil2014stabilised}\cite{AGUIRRE2015387}.
In recent works \cite{bonet2015first}\cite{gil2016first}\cite{HASSAN2019100025}, a nodal co-factor tensor $\textbf{H}= cof  :\ \textbf{F}$ is added.
This method, like others in the family of methods based on nodal interpolations of $\textbf{F}$ are inherently unstable. 
Thus, a stabilization based on the Streamline Upwind/Petrov-Galerkin (SUPG) method, and added penalties on the deformation gradient $\textbf{F}$ is utilized.
Tests show the capabilities of this methodology to solve problems in the incompressible limit, it is however costly in terms of the number of unknowns per node.
\\
\\
\indent
In \cite{articleSC}\cite{articleRO}, a mixed problem in velocity and pressure, where the displacement field is calculated based on the discretization of the velocity is presented.
Authors claim that the Variational Multi-Scale (VMS) method was insufficient for the hyperbolic problem of the transient case, so they resorted to a pressure rate equation to elevate the problem.
The tests also prove to be accurate and robust.
In \cite{articleCA}, a mixed formulation in displacement and pressure, resolved in the total Lagrangian framework was presented with different variations of the VMS methods.
The method proved to be accurate and robust as well.
\\
\\
\indent
In \cite{MASUD20113453}, a mixed problem in displacement and  pressure in finite elements for nearly incompressible material is presented.
The Variational Multi-Scale (VMS) is used for the displacement field, and two types of error estimators are exploited. 
The formulation was investigated across different numerical convergence tests.
\\
\\
\indent
This paper proposes a novel method for solving non-linear elasticity in solid dynamics. 
Based on unstructured tetrahedral meshes, the method is able to depict complex geometries with ease and acceptable computational cost.
In addition to the split of the strain energy into its deviatoric and volumetric part \cite{SANSOUR200828}, a constitutive equation in pressure is also solved.
A fully implicit, mixed coupled in displacement and pressure (piece-wise linear) formulation in an updated Lagrangian context is proposed.
The set of equations obtained is prone to spurious pressure fields.
A stabilization based on the Variational Multi-Scale (VMS) method is thus implemented to elevate the problem.
This method can handle complex geometries with a reasonable computational time.
\\
\\
\indent
The rest of the paper consists of: Section 2 contains the problem definition; Section 3 presents the stabilized linear elastic formulation both in its steady state and transient form; Section 4 contains the stabilized hyperelastic formulation; Section 5 provides the numerical validation of the framework.
Finally, perspective and conclusions are given in section 6.

\section{Lagrangian Solid Dynamics}
\indent

\subsection{Solid Dynamics}

The variation rate of density, and displacement for a solid material is governed by the equations of Lagrangian solid dynamics.
$\Omega_0$, and $\Omega$ represent the initial and current domain, which are two open sets in ${\rm I\!R}^d$ with Lipshitz boundaries, where $d$ denotes the spatial dimension.
The boundary is given by $\Gamma$, which is split into two separate sets given by $\Gamma=\overline{\partial\Omega_{\uu} \cup \partial\Omega_{\textbf{t}}}$ and $\partial\Omega_{\uu} \cap \partial \Omega_{\textbf{t}} = \emptyset$.
Where $\partial\Omega_{\uu}$ denotes the dirichlet boundary that specifies the displacement, and  $\partial\Omega_{\textbf{t}}$ denotes the Neumann boundary that specifies the traction force.
The motion of the deformable body is given by:

\begin{equation}
\boldsymbol{\phi} : = \Omega_0 \rightarrow \Omega = \boldsymbol{\phi} (\Omega_0) 
\end{equation}

\begin{equation}
\boldsymbol{\phi} : = \Gamma_0 \rightarrow \Gamma = \boldsymbol{\phi} (\Gamma_0) 
\end{equation}

\begin{equation}
\mathbf{X} \rightarrow \mathbf{x} = \boldsymbol{\phi}(\mathbf{X},t) \: \forall \mathbf{X} \in \Omega_0
\end{equation}

It serves as a mapping of the material coordinate $\mathbf{X}$, in the total Lagrangian framework of an infinitesimal material particle of the solid, to $ \mathbf{x}$, the coordinate of the same particle in the updated Lagrangian framework.
$\boldsymbol{\phi}$ is assumed to be smooth, and invertible.
The deformation gradient and the Jacobian determinant are given by: $ \mathbf{F} =  \nabla_{\mathbf{X}} \boldsymbol{\phi}$ and $ J = det \mathbf{F}$.
\\
\\
\indent
The displacement of the solid is given by: $ \mathbf{u}= \mathbf{x} - \mathbf{X}$. 
The governing equations are given by:

\begin{equation}
\rho \ddot{\uu} = \nabla . \: \bfsigma + \bff \: in \: \Omega
\end{equation}
\begin{equation}
\rho J = \rho_0 \: on \: \partial \Omega_u
\end{equation}

Where $\rho$, and  $\rho_0$ are the current and initial body density respectively, $ \bff $ is a forcing term,  $\bfsigma$ is the symmetric Cauchy stress tensor, and the derivatives are taken with the respect to updated reference frame.
Moreover, $ \ddot{\uu}$ represents the material second derivative of displacement, which is the acceleration.
This set of equations, along with a constitutive model for the solid that defines $\bfsigma$, and the corresponding initial and boundary condition, describes the development of the system.
\\
\\
\indent
One way to model a solid is using a mixed formulation, containing both displacement and pressure fields.
This is obtained with a decomposition of the stress into a volumetric and deviatoric component.
This decomposition is essential when dealing with incompressible or nearly incompressible material. 
It is significant to mention that this is done for isotropic material.
Thus the stress is given by:

\begin{equation}
\bfsigma = p \textbf{I} + dev[\bfsigma] 
\end{equation}

Where $\textbf{I}$ is the identity matrix
\\
\\
\indent
The problem is completed with the addition of the initial and boundary conditions of the problem.
Assuming zero displacement initial conditions, given by $ \uu(\mathbb{X},0)=\uu_0 = 0$.
This gives: $\boldsymbol{\phi}(\mathbf{X},0) = \mathbf{X}$, $\textbf{F}|_{t=0} = \textbf{I}$, and $J|_{t=0}=1$.
The material is also assumed to be stress free.
The boundary condition are given by:

\begin{equation}
\uu|_{\Gamma_{\uu}}= \uu(\xx,t)
\end{equation}
\begin{equation}
\bfsigma\textbf{n}|_{\Gamma_{\textbf{t}}}= \textbf{t}(\xx,t)
\end{equation}

Where $\textbf{n}$ is the outward-pointing normal on the boundary $\Gamma$.

\section{Linear Elastic Formulation}
\indent

For very small displacement, the elastic behavior of the solid can be modeled using Hook's law.
By considering a linear relationship between stress and strain, the solid is modeled using a spring.
A steady-State formulation is first introduced for comparison purposes, and a transient formulation is afterwards developed. 

\subsection{Steady-State Formulation}

To put the different implementations that are going to be presented in this paper into perspective, the steady-state linear elastic solver is first developed.
The latter is a direct extension of the stokes problem, where the displacement $\uu$ is the primary variable instead of the velocity.
This extension will better explain the stabilization techniques that are being exploited.
For linear elasticity, $\mathbf{x} \approx \mathbf{X}$, $\Omega_0 \approx \Omega$,
$\nabla_\mathbb{X} \approx \nabla_\mathbf{x} $, and $\rho_0=\rho$.
A linear elastic problem can be formulated based on the decomposition of the stress, using the hydro-static pressure $p$ and the displacement field $\uu$.
It is a worthwhile to note that the pressure convention in solid mechanics is opposite to that of fluid mechanics.

The stress tensor is thus given by:
\begin{equation}
\label{eqn:9}
\bfsigma = p \textbf{I} +2\mu \: dev[\nabla^s \uu] 
\end{equation}
\begin{equation}
\label{eqn:10}
p=K \epsilon_v
\end{equation}
\begin{equation}
\label{eqn:11}
\epsilon_v = \nabla . \uu
\end{equation}

In equations (\ref{eqn:9})(\ref{eqn:10})(\ref{eqn:11}), we distinguish between the deviatoric and volumetric part of the deformation.
$\nabla^s$ is the symmetrical gradient operator:
\begin{equation}
\nabla^s = \frac{1}{2}(\nabla + \nabla^T)
\end{equation}

$\mu$ is the Lam\'e constant, also know as the shear modulus of the material, and it is specified by:
\begin{equation}
\mu=\frac{E}{2(1+\nu)}
\end{equation}
K is the bulk modulus or modulus of volumetric compressibility, and it is defined by:
\begin{equation}
K=\frac{1}{3}\frac{E}{(1-2\nu)}
\end{equation}

Using the stress tensor formulation, along with a body force $\bff$ and the necessary dirichlet and neumann boundary conditions, the steady state problem can be formulated as follows:
\begin{equation}
\nabla p + 2\mu \nabla . \: dev[\nabla^s \uu] + \bff = 0 \: in \: \Omega
\end{equation}
\begin{equation}
\frac{1}{K} p - \nabla . \uu=0 \: in \: \Omega
\end{equation}
\begin{equation}
\uu=0 \: on \: \partial \Omega_u
\end{equation}
\begin{equation}
\bfsigma.\bfn=\overline{\bft} \: on \: \partial \Omega_t
\end{equation}
\begin{equation}
\rho \mathbf{J} = \rho_0 \: in \: \Omega
\end{equation}

This formulation takes into account both incompressible and compressible material, with the difference being in equation (16) and a constant density $\rho_0$.
For an incompressible material, $K \rightarrow \infty$ and equation (16) becomes simply:
\begin{equation}
\nabla . \uu = 0 \: in \: \Omega
\end{equation}

Whereas if we assume an isochoric phenomena, implying $\epsilon_v = 0$, we will get the same result.

The variational formulation of this problem is given by:

\begin{equation}
a(\uu,\vv)+(p,\nabla . \vv)=L(\vv) \: \forall \vv \: \in V_0
\end{equation}
\begin{equation}
(\nabla . \uu, q)-(\frac{1}{K} p,q)=0 \: \forall q \: \in Q
\end{equation}

Where $a(\uu,\vv)$ and $L(\vv)$ are given by:

\begin{equation}
a(\uu,\vv)=\int_\Omega 2\mu \: dev[\nabla^s \uu] : \nabla^s \vv \: d\Omega
\end{equation}
\begin{equation}
L(\vv)=\int_\Omega \bff . \vv d\Omega + \int_{\partial \Omega_t} \vv . \overline{\bft} d\Gamma
\end{equation}

The discrete form of the problem is given by:

\begin{equation}
a(\uu_h,\vv_h)+(p_h,\nabla . \vv_h)=L(\vv_h) \: \forall \vv_h \: \in V_{h,0}
\end{equation}

\begin{equation}
(\nabla . \uu_h, q_h)-(\frac{1}{K} p_h,q_h)=0 \: \forall q_h \: \in Q_h
\end{equation}

The Babuska--Brezzi or inf-sup stability \cite{babuvska1971error} condition constrains the interpolation relation between the fields, thus forcing different interpolations for $\uu$ and $p$.
Equal order interpolation has poor numerical performance as it does not respect the condition.
Several types of stabilization are available in this case.
We used P1/P1 elements, with a Variational MultiScale Method (VMS), which enable us to have the same order of interpolation.
In \cite{HUGHES198685}, equal order elements were used for the Stokes problem. 
It contained proof of convergence, and stability
This work led to the extension of the formulation to the Navier--Stokes equations \cite{HACHEM20108643}.
In \cite{HUGHES198785} \cite{10.1007/BF01395881}, the linear elastic problem was tackled.
VMS provides natural stabilization by an orthogonal decomposition of the solution (displacement, pressure) spaces.
Orthogonal decomposition of the function spaces is first done by:

\begin{equation}
V_0 = V_{h,0} + V_0'
\end{equation}
\begin{equation}
Q = Q_h + Q'
\end{equation}

Following \cite{hughes1998variational}, the resolvable coarse and unresolved fine scale components of the displacement and pressure are given by:

\begin{equation}
\uu=\uu_h+\uu'
\end{equation}
\begin{equation}
p=p_h+p'
\end{equation}

We also apply the same decomposition for the weighting functions:

\begin{equation}
\vv=\vv_h+\vv'
\end{equation}
\begin{equation}
q=q_h+q'
\end{equation}

The equations are divided into two sets: coarse, and fine scale.
The unresolved fine-scales, are most of the time modeled in function of the residual based terms.
Using static condensation, the fine scale equations are solved in an approximate matter (residual based) and re-injected into the coarse scale equations.
This will provide us with additional terms, calibrated by a local stabilizing parameter.
These terms are responsible for the enhanced stability, reduced pressure oscillations, and increased accuracy of the standard Galerkin formulation.

The fine scale problem, defined on the sum of elements interiors \cite{BREZZI1997329}, and formulated in function of the transient coarse scale variables, is solved.
The fine scale approximations are given by:

\begin{equation}
\uu' = \sum_{T_h} (\tau_{\uu} P_{\uu}' (R_{\uu})
\end{equation}
\begin{equation}
p' = \sum_{T_c} (\tau_c P_{c}' (R_c)
\end{equation}

Where $R_{\uu}$, and $R_c$ are the finite elements residuals, $P_{\uu}'$, and $P_{c}'$ are the projection operators, and $\tau_{\uu}$, and $\tau_c$ are the tuning parameters.
Note that in this current work, both $P_{\uu}'$, and $P_{c}'$ are taken as the Identity matrix.

The fine scale approximations are subsequently substituted in the coarse problem.
The new variational formulation for the coarse scale equations are given by:

\begin{equation}
a((\uu_h+\uu'),\vv_h)+(p_h+p', \nabla.\vv_h)=L(\vv_h) \: \forall \: \vv_h \: \in \: V_{h,0}
\end{equation}
\begin{equation}
(\nabla.(\uu_h+\uu'),q_h)-(\frac{1}{K}(p_h+p'),q_h)=0  \: \forall \: q_h \: \in \: Q_{h}
\end{equation}

and that of the fine scale equations are given by:

\begin{equation}
a((\uu_h+\uu'),\vv')+(p', \nabla.\vv')=L(\vv') \: \forall \: \vv' \: \in \: V_{0}'
\end{equation}
\begin{equation}
(\nabla.(\uu_h+\uu'),q')-(\frac{1}{K}(p_h),q')=0  \: \forall \: q' \: \in \: Q'
\end{equation}

Finally, calculating the fine scale equation based on the initial residual, and re-entering the physics in the coarse scale equation, we get the final set of the coarse scale equations with the pressure stabilization term for the case of linear elasticity, given by:

\begin{equation}
a((\uu_h, \vv_h)+(p_h, \nabla .\vv_h) - L(\vv_h)=0 \: \forall \: \vv_h \in V_h
\end{equation}
\begin{equation}
(\nabla . \uu_h, q_h) - (\frac{1}{K} p_h, q_h) + \sum_{K \in T_h} (\tau_K R (\uu_h),\nabla q_h) \: \forall \: q_h \in Q_h
\end{equation} 
\begin{equation}
R (\uu_h) = \nabla p_h + 2\mu \nabla . \: dev[\nabla^s \uu_h] + \bff
\end{equation} 

Where $R (\uu_h)$ is the finite element residual, and $\tau_K$ is a coefficient based on the study of the response of the stabilization parameters coming from a Fourier analysis of the problem for the sub-scales \cite{cervera2010mixed}.

Comparing the standard Galerkin and the stabilized formulation, we distinguish additional integrals that are evaluated element-wise.
These terms represent the sub-grid scales, and help damp out spurious pressure oscillations, and overcome instabilities in our case.

\subsection{Transient Formulation}

A second order derivative in time for the displacement, which accounts for the dynamics of the solid is added.
While the previous steady-state equations are parabolic in nature, the added transient term renders the PDE hyperbolic.
When dealing with materials in the incompressible limit, the PDE becomes degenerate hyperbolic.
The pressure acts as a Lagrangian multiplier required to force the divergence-free constraint of the displacement.
The transient elastic solid solver governing equations are given by:

\begin{equation}
\rho \ddot{\uu} = \nabla p + 2\mu \nabla . \: dev[\nabla^s \uu] + \bff \: in \: \Omega
\end{equation}
\begin{equation}
\frac{1}{K} p - \nabla . \uu=0 \: in \: \Omega
\end{equation}
\begin{equation}
\uu=0 \: on \: \partial \Omega_u
\end{equation}
\begin{equation}
\bfsigma.\bfn=\overline{\bft} \: on \: \partial \Omega_t
\end{equation}

The Courant-Friedrichs-Lewy (CFL) condition, imposes limits on the time step.
For explicit time integrators, very small time steps are needed to obtain accurate results. 
Consider a time interval, where $t \in [0,T]$, and a discretization of this interval into $N$ time steps ($\Delta t$).
A Backward differentiation formula (BDF) is adopted in this work.
A first, and second order accurate BDF's are considered.
These equations are given by:

\begin{equation}
\textbf{a}^{n+1} \approx \frac{1}{\Delta t^2} ( \uu^{n+1} -2 \uu^{n} + \uu^{n-1}) + O(\Delta t)
\end{equation}

\begin{equation}
\textbf{a}^{n+1} \approx \frac{1}{\Delta t^2} ( 2\uu^{n+1} -5 \uu^{n} + 4\uu^{n-1} - \uu^{n-2}) + O(\Delta t^2)
\end{equation}

The backward differentiation formulas are known for their high frequency dissipation, that will help damp out spurious high frequency oscillations.

The discrete/stabilized form in the variational form of the above equation, following the same steps as before, are given by:

\begin{equation}
(\rho \ddot{\uu},\vv_h) + a(\uu_h, \vv_h)+(p_h, \nabla .\vv_h) - L(\vv_h)=0 \: \forall \: \vv_h \in V_h
\end{equation}
\begin{equation}
(\nabla . \uu_h, q_h) - (\frac{1}{K} p_h, q_h) + \sum_{K \in T_h} (\tau_K R (\uu_h),\nabla q_h)=0 \: \forall \: q_h \in Q_h
\end{equation} 
\begin{equation}
R (\uu_h) = - \rho \ddot{\uu} +  \nabla p_h + 2\mu \nabla . \: dev[\nabla^s \uu_h] + \bff
\end{equation}

Where $R (\uu_h)$ is the new finite element residual.
The same additional elements that damp out pressure oscillations are found in these equations as well.

\section{Transient Non-Linear Elastic Lagrangian Formulation}
\subsection{Hyperelasticity model and pressure equation}

Elastic materials in general are better modeled with a non-linear depiction of their real life behavior.
As mentioned earlier, the material at hand is considered to be isotropic.
Consider a nonlinear material with a Helmholtz free energy or strain energy $\Psi (\textbf{C})$ function, where $\textbf{C}$ is the right Cauchy-Green strain tensor $\textbf{C} = \textbf{F}^T\textbf{F}$.
$\textbf{F}$ is the deformation gradient given by $\textbf{F}_{ij} = x_{i,j} = \frac{\partial x_i}{\partial X_j}$.
The second Piola-Kirchhoff stress tensor $\textbf{S} =  \textbf{J} \textbf{F}^{-1} \bfsigma \textbf{F}^{-T}$ where $\textbf{J}$ is the Jacobian determinant of F, is derived by taking derivatives of the Helmholtz free energy functional  $\boldsymbol{\Psi} (\textbf{C})$ with respect to $\textbf{C}$:

\begin{equation}
\textbf{S}=2\partial_{\textbf{C}} \Psi (\textbf{C}).
\end{equation}

To be able to model both incompressible and compressible material, we will apply the same decomposition as before.
Decomposing $\Psi ( \textbf{C} )$ into its volumetric and deviatoric part respectively as follows:

\begin{equation}
\Psi ( \textbf{C} ) = U(\textbf{J}) + W( \bar{\textbf{C}} )
\end{equation}

Where $\bar{\textbf{C}} = \textbf{J}^{-\frac{2}{3}} \textbf{C}$ is the deviatoric/volume-preserving part of $\textbf{C}$, and $\textbf{J} = \sqrt{det \textbf{C}}$.

The Helmholtz free energy of isotropic hyperelastic models is written in function of the strain invariants.
Consider a Neo-Hookean elastic material, and a Simo-Taylor volumetric model \cite{simo1985variational} with:

\begin{equation}
U( \textbf{J} ) = \frac{1}{4} \kappa ( \textbf{J}^2 - 1 ) - \frac{1}{2} \kappa ln \textbf{J}
\end{equation} 

\begin{equation}
W ( \bar{\textbf{C}} ) = \frac{1}{2} \mu ( tr\bar{\textbf{C}} - 3 ) 
\end{equation}

where $\kappa$ and $\mu$ are material properties. 
For small displacements, the model reduces to a linear elastic model where $\kappa$ and $\mu$ are the bulk and shear modulus of material.
The stress can also be split to its deviatoric and volumetric part:

\begin{equation}
p = 2 \textbf{J}^{-1} \textbf{F} \frac{\partial U(\textbf{J})}{\partial \textbf{C}} \textbf{F}^T = U'(\textbf{J}) = \frac{1}{2} \kappa (\textbf{J} + \textbf{J}^{-1})
\end{equation}

\begin{equation}
\textbf{dev}[\bfsigma] = 2 \textbf{J}^{-1} \textbf{F} \frac{\partial W(\bar{\textbf{C}})}{\partial \textbf{C}} \textbf{F}^T = \mu \textbf{J}^{-\frac{5}{3}} \textbf{dev} [ \textbf{FF}^T]
\end{equation}

Recall, $ \mathbf{F} =  \nabla_{\mathbf{X}} \uu + I$. Thus:

\begin{equation}
\mathbf{F} \mathbf{F^T} =  \nabla_{\mathbf{X}} \uu +  \nabla_{\mathbf{X}}^T \uu + \nabla_{\mathbf{X}} \uu  \nabla_{\mathbf{X}}^T \uu + I
\end{equation}

We are solving our equations in the updated Lagrangian framework, while the above equation is given in the total Lagrangian framework.
Considering the following mathematical equation:

\begin{equation}
\nabla_{\mathbf{X}} \uu = ( I - \nabla \uu ) ^{-1} - I
\end{equation}

And assuming a very small variation in the displacement noted $\delta \uu$.
Recalling that for very small displacement, $( I - \nabla \uu ) ^{-1} = I + \nabla \uu$, we get:

\begin{equation}
\begin{split}
\mathbf{F} \mathbf{F^T} =  & ( I - \nabla \uu ) ^{-1} - I + ( ( I - \nabla \uu ) ^{-1} - I)^T +  (( I - \nabla \uu ) ^{-1} - I)( ( I - \nabla \uu ) ^{-1} - I)^T  \\
& + I + 2\epsilon (\delta \uu) + \nabla \delta \uu (\nabla \delta \uu)^T + \nabla \delta \uu (\nabla  \uu)^T + (\nabla  \uu) ( \nabla \delta \uu)^T
\end{split}
\end{equation}

The system of equations to be solved now is given by:

\begin{equation}
\rho \ddot{\uu} = \nabla_x p +  \nabla_x . \: dev[\bfsigma] + \bff \: in \: \Omega
\end{equation}

\begin{equation}
\frac{1}{K} p -  \nabla_x . \uu=0 \: in \: \Omega
\end{equation}

\begin{equation}
\uu=0 \: on \: \partial \Omega_u
\end{equation}

\begin{equation}
\bfsigma.\bfn=\overline{\bft} \: on \: \partial \Omega_t
\end{equation}

\begin{equation}
\rho J = \rho_0
\end{equation}

$\ddot{\uu}$ represents the time material derivative of displacement.
The formulation becomes:

\begin{equation}
\rho  \Bigg(\frac{\partial^2 \uu}{\partial t^2} + (\nabla_x( \frac{\uu - \uu^-}{\Delta t}).(\vv - \vv_{domain})\Bigg)= \nabla_x p +  \nabla_x . \: dev[\bfsigma] + \bff \: in \: \Omega
\end{equation}

\section{Moving Mesh Method (MMM)}

For $\vv_{domain}$, we adopt the R-method as an adaptive strategy \cite{tang2005moving}.
The r-method or moving mesh method (MMM), consists of relocating mesh nodes so that nodes get condensed in regions with high gradients.
This is done, through a mapping from the undeformed domain in a parameter space $\Omega_c$, to the deformed domain in the physical space $\Omega$.
The connections of points in $\Omega$, representing discrete points in $\Omega_c$, ensures the coverage of the physical domain with a computational mesh.
The key components are threefold: 

\begin{enumerate}

\item Mesh equations
\item Monitor Function
\item Interpolation

\end{enumerate}

Choosing the appropriate mesh equations for a given application and resolving them efficiently is essential for the method.
In our case the mesh is guided by the solid dynamics equation.
Guiding the mesh redistribution is done via the monitor function.
It is dependent on the solutions arclength in 1D, curvature, and a posteriori errors if needed.
It also requires smoothing in practice.
Interpolation of dependent variables from the old, to the new mesh is only needed if the mesh equations are not time dependent and are solved independently from the partial differential equation. 

In interpolation free MMM, such as the moving finite element method of Miller \cite{10.2307/2157255}
\cite{10.2307/2157254}, both the differential equation and the mesh equations are resolved simultaneously.
The essential components of such methods include:

\begin{enumerate}

\item Equidistribution principle
\item Mesh equations
\item The method of lines (MOL) approach

\end{enumerate}

First introduced by de Boor in \cite{10.1007/BFb0069121}, the equidistribution principle  was used to solve Boundary Value Problems (BVP) for Ordinary Differential Equations (ODE).
It consists of choosing mesh nodes, so that a certain measure describing the solution error is adjusted over every sub-interval.
The method of lines approach is usually considered in most moving mesh codes, which may result in a stiff equation.
A moving mesh finite element approach is used in our case, where the mesh equations are based on the solid dynamics equation, and the equidistribution principle is applied on the residual of the partial differential equation written in finite element form.
This method is particularly interesting for its various advantages:

\begin{itemize}

\item Interpolation free method.
\item Detect, track and resolve moving boundaries.
\item The method of lines (MOL) approach

\end{itemize}

\begin{figure}[H]
\centering
   \includegraphics[width=0.7\linewidth]{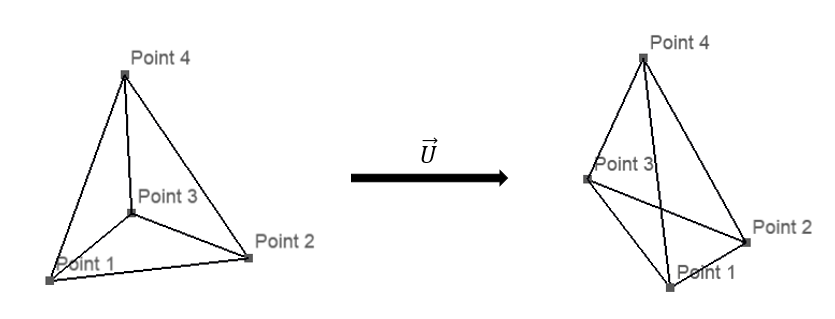}
\caption{Moving mesh illustration with varying volume in 3D}
\label{fig:MMM}
\end{figure}

It is however necessary to use an implicit time scheme to overcome the stiffness of the system.
An illustration of an element subjected to a displacement vector is shown in figure \ref{fig:MMM}.

Since this method is adpoted, thus we move the mesh with the velocity of the solid.
$\vv=\vv_{domain} \rightarrow \vv-\vv_{domain}= 0$.
The variational formulations thus become:

\begin{equation}
(\rho \frac{\partial^2 \uu}{\partial t^2},\ww_h) + a'((\uu_h, \ww_h)+(p_h, \nabla .\ww_h) - L(\ww_h)=0 \: \forall \: \ww_h \in W_h
\end{equation}
\begin{equation}
(\nabla . \uu_h, q_h) - (\frac{1}{K} p_h, q_h)=0 \: \forall \: q_h \in Q_h
\end{equation} 

where $a'$ is given by:

\begin{equation}
a'((\uu_h, \ww_h)= \int_\Omega \mu \: dev[\boldsymbol{\sigma}] : \nabla^s \ww \: d\Omega
\end{equation}

Applying the VMS method for this equation, we end up with:

\begin{equation}
\begin{split}
(\rho \frac{\partial^2 \uu}{\partial t^2},\ww_h) + a'((\uu_h, \ww_h) - (p_h, \nabla .\ww_h) - L(\ww_h ) + (\nabla . \uu_h, q_h) - (\frac{1}{K} p_h, q_h) \\ + \sum_{K \in T_h} (\tau_K R (\uu_h),\nabla q_h) + \sum_{K \in T_h} (\tau_K \nabla . (\uu_h),\nabla . (\ww_h)) - \sum_{K \in T_h} (\frac{\tau_K}{K} p_h,\nabla . (\ww_h)) = 0 
\end{split}
\end{equation}

Where $R (\uu_h)$ is the new finite element residual.

This is a more general stabilization formulation than that of the linear elastic case, as it contains two additional terms, that help impose the incompressibility constraint.

\section{Numerical Validation}
\subsection{Linear Elastic}
\subparagraph{Static Cook's membrane test}

\begin{figure}
\centering
\begin{subfigure}{0.4\textwidth}
   \includegraphics[width=1\linewidth]{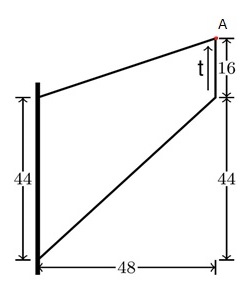}
    \caption{Geometrical Steup} \label{fig:1a}
  \end{subfigure}%
  \hspace*{\fill}   % maximize separation between the subfigures
  \begin{subfigure}{0.5\textwidth}
   \includegraphics[width=1\linewidth]{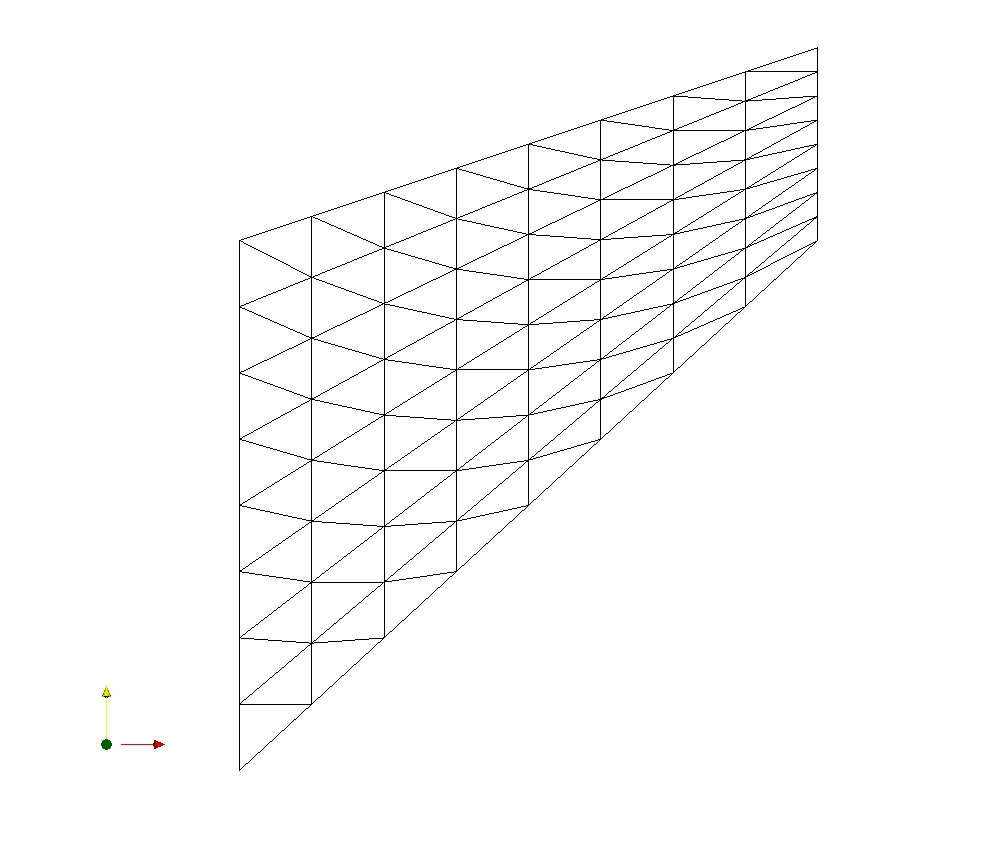}
    \caption{Structured mesh} \label{fig:1b}
  \end{subfigure}%
  \\
  \centering
  \begin{subfigure}{0.5\textwidth}
   \includegraphics[width=1\linewidth]{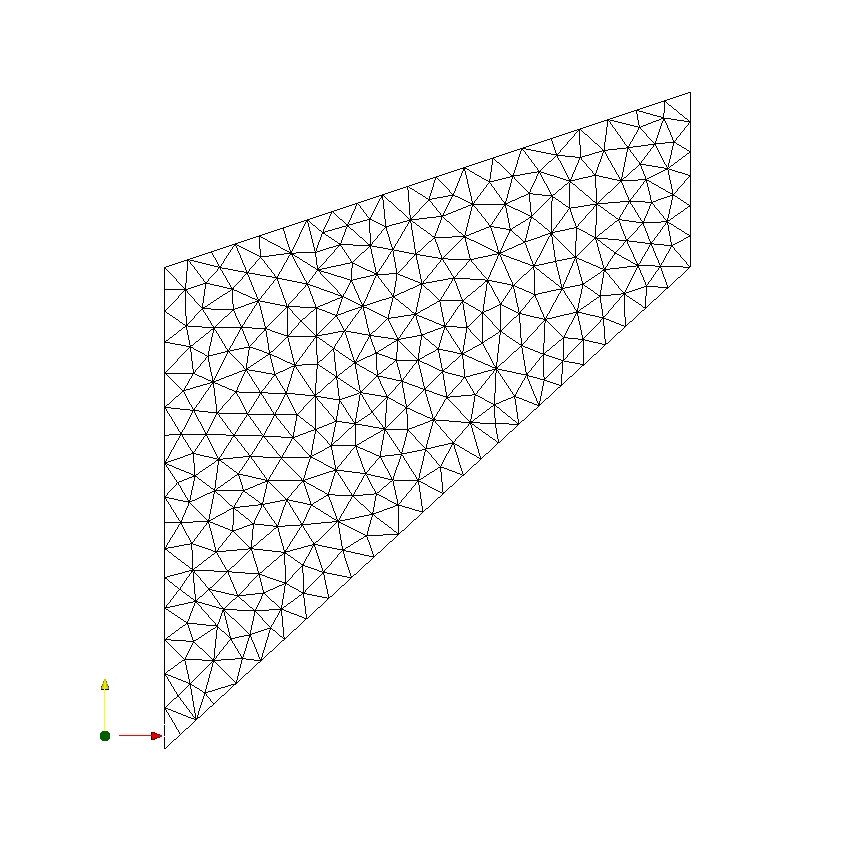}
    \caption{Unstructured mesh} \label{fig:1c}
  \end{subfigure}

\caption{Geometrical Setup and 2D meshs}
\label{fig:GeoMesh}
\end{figure}

A typical problem where the P1/P1 elements for both displacement and pressure produce a polluted pressure field is the Cook's membrane problem \cite{cook1974improved}.
The geometrical setup of the problem and a typical structured, and unstructured 2D mesh are shown in figure \ref{fig:GeoMesh}.

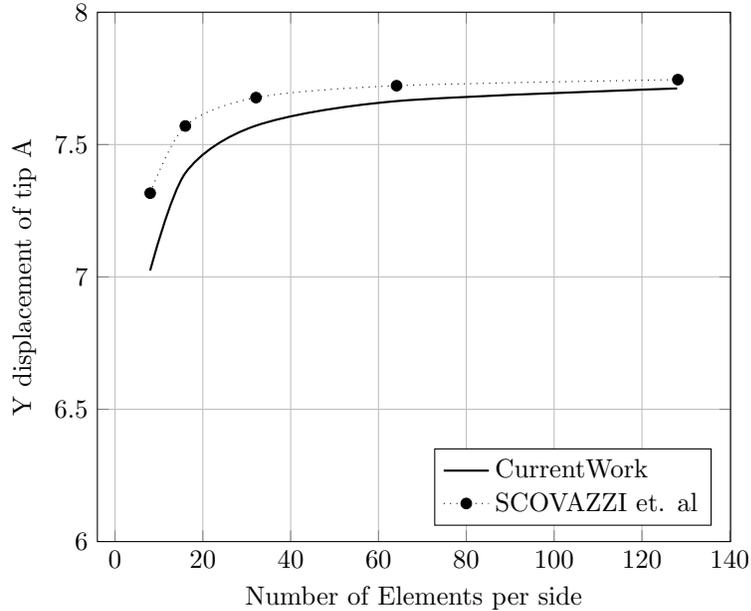
\begin{figure}
        \centering
        \begin{tikzpicture}[trim axis left, trim axis right]
                \begin{axis}[
                legend cell align=left,
                legend pos=south east,
                grid=major,
                ymin=6,
                ymax=8,
                xlabel={Number of Elements per side},
                ylabel={Y displacement of tip A}
                ]                 
                 
                \legend{CurrentWork, SCOVAZZI et. al}
      
                \addplot[ 
                draw=black,
                thick,
                smooth
                ]
                table[
                x index=0,
                y index=1
                ]
                {data/Ydisp.dat};

                \addplot[ 
                color=black,
                dotted,
                mark=*,
                mark options={solid},
                smooth
                ]
                 table[
                x index=0,
                y index=1
                ]
                {data/YdispBenchmark.dat};  
               
                \end{axis}
        \end{tikzpicture}
        \caption{Mesh convergence study}
        \label{fig:MeshConvergence}
\end{figure}

\begin{figure}
\centering
\begin{subfigure}{0.4\textwidth}
   \includegraphics[width=1\linewidth]{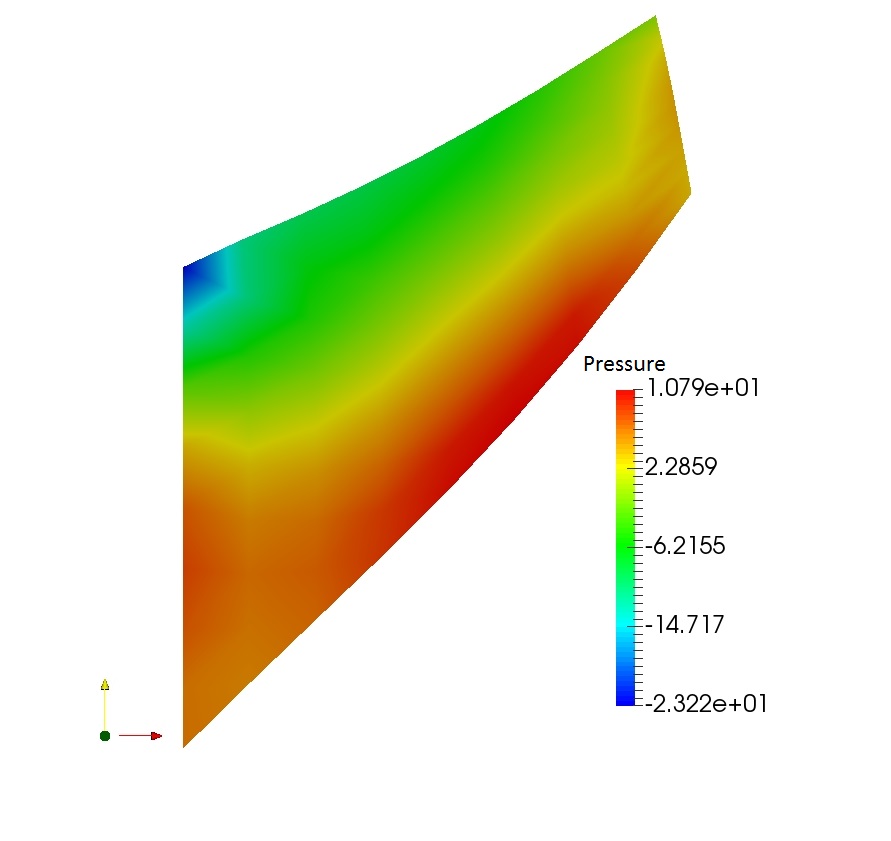}
    \caption{2x8x8 structured mesh pressure contours} \label{fig:3a}
  \end{subfigure}%
  \hspace*{\fill}   % maximize separation between the subfigures
  \begin{subfigure}{0.4\textwidth}
   \includegraphics[width=1\linewidth]{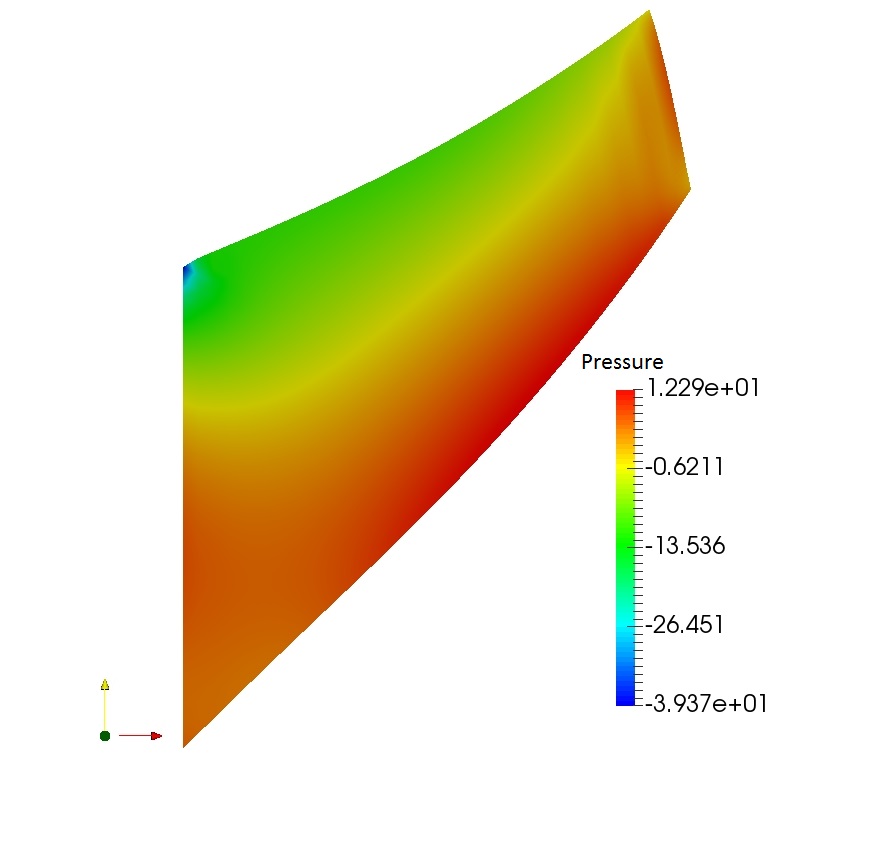}
    \caption{2x32x32 structured mesh pressure contours} \label{fig:3b}
  \end{subfigure}%
    \\
  \centering
  \begin{subfigure}{0.4\textwidth}
   \includegraphics[width=1\linewidth]{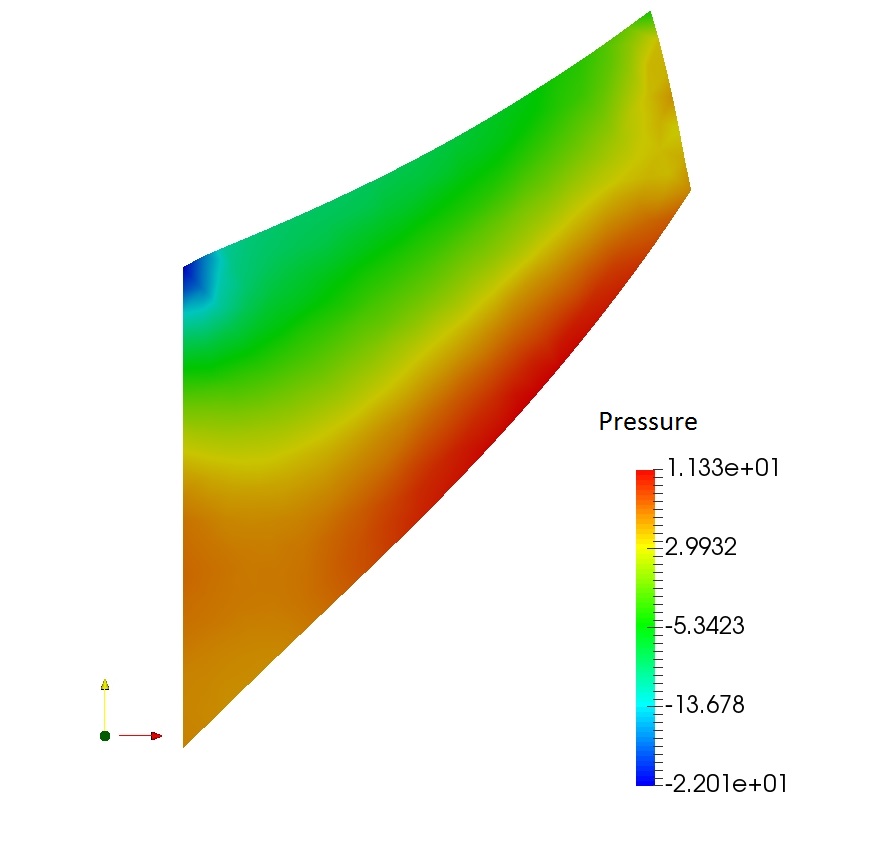}
    \caption{Unstructured mesh pressure contours} \label{fig:3c}
  \end{subfigure}

\caption{Pressure Contours for different meshes}
\label{fig:PressureContours}
\end{figure}

The material is assumed to have a linear elastic behavior with $\rho=1$, $E=250$, and $\nu=0.49995$. Zero displacement dirichlet boundary conditions are imposed on the left side of the membrane, and a uniform vertical traction force equal to 6.25 is imposed on the right side of the membrane.
Standard Galerkin P1/P1 elements for both the displacement and pressure lead to oscillations in pressure when no stabilization is included. 
However, when the VMS stabilization is applied the pressure field has no oscillatory pulses in the pressure.
Different pressure contours are found in figure [\ref{fig:PressureContours}].
It is important to note that we get similar result when using unstructured tetrahedral meshes.
A mesh convergence study was applied on refiened structured meshes, and the results obtained conform with the literature and can be found in figure [\ref{fig:MeshConvergence}].

\subparagraph{Transient Cook's membrane test}

\begin{figure}
\centering
\begin{subfigure}{0.48\textwidth}
   \includegraphics[width=1\linewidth]{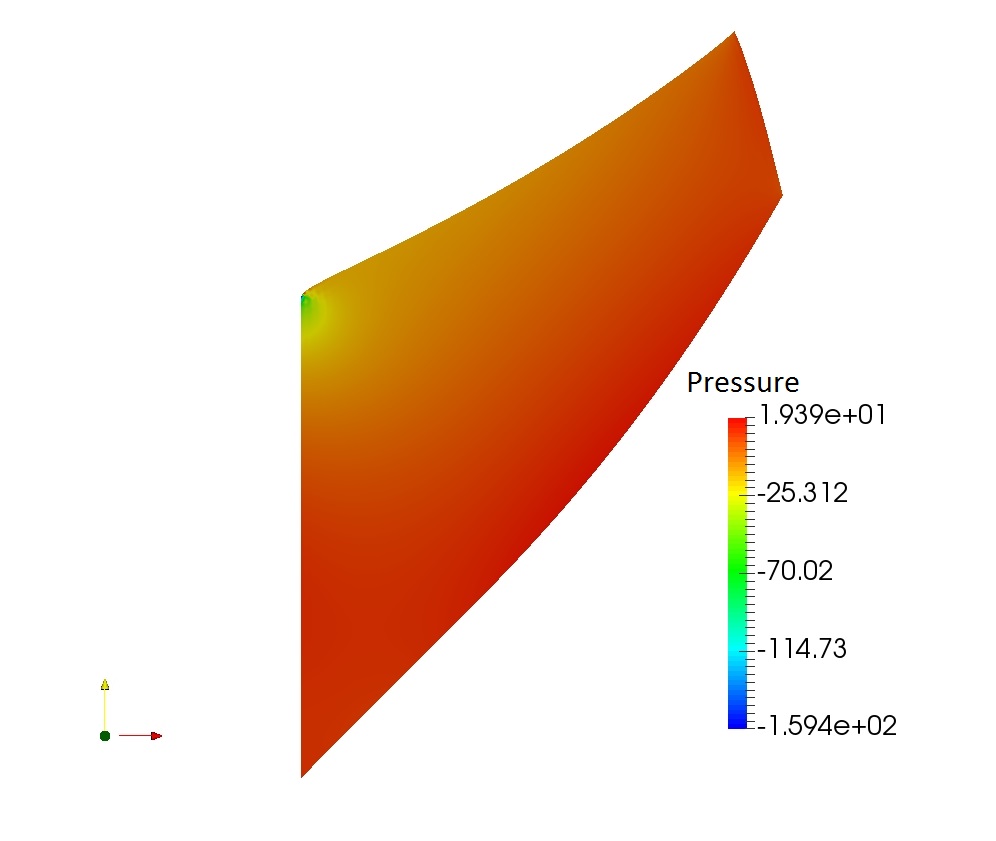}
    \caption{t=1s} \label{fig:2a}
  \end{subfigure}%
  \hspace*{\fill}   % maximize separation between the subfigures
  \begin{subfigure}{0.48\textwidth}
   \includegraphics[width=1\linewidth]{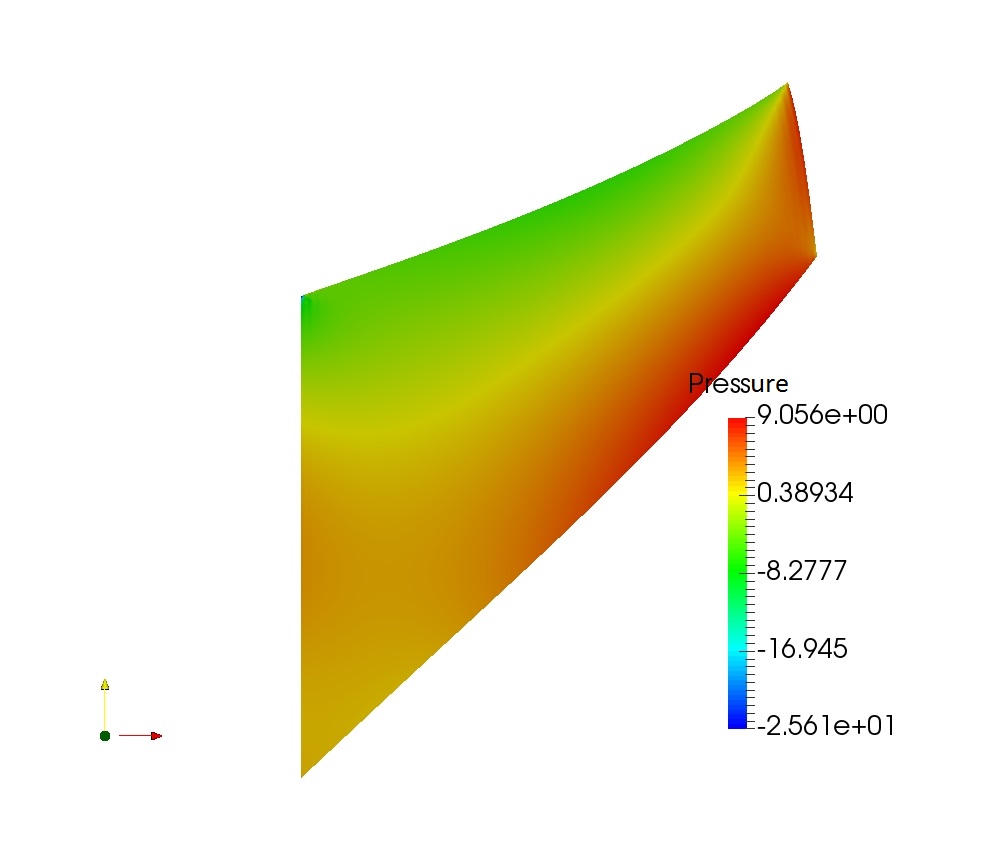}
    \caption{t=3s} \label{fig:2b}
  \end{subfigure}%
  \\
  \hspace*{\fill}   % maximize separation between the subfigures
  \begin{subfigure}{0.48\textwidth}
   \includegraphics[width=1\linewidth]{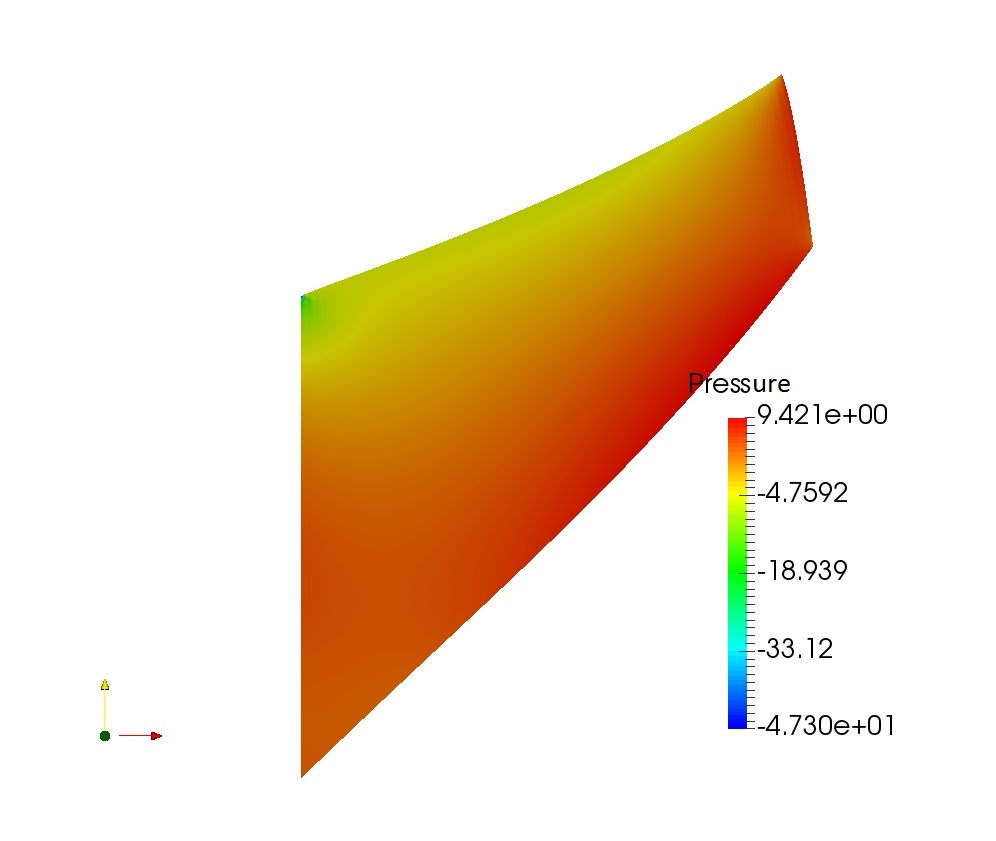}
    \caption{t=5s} \label{fig:2c}
  \end{subfigure}
  \hspace*{\fill}   % maximize separation between the subfigures
  \begin{subfigure}{0.48\textwidth}
   \includegraphics[width=1\linewidth]{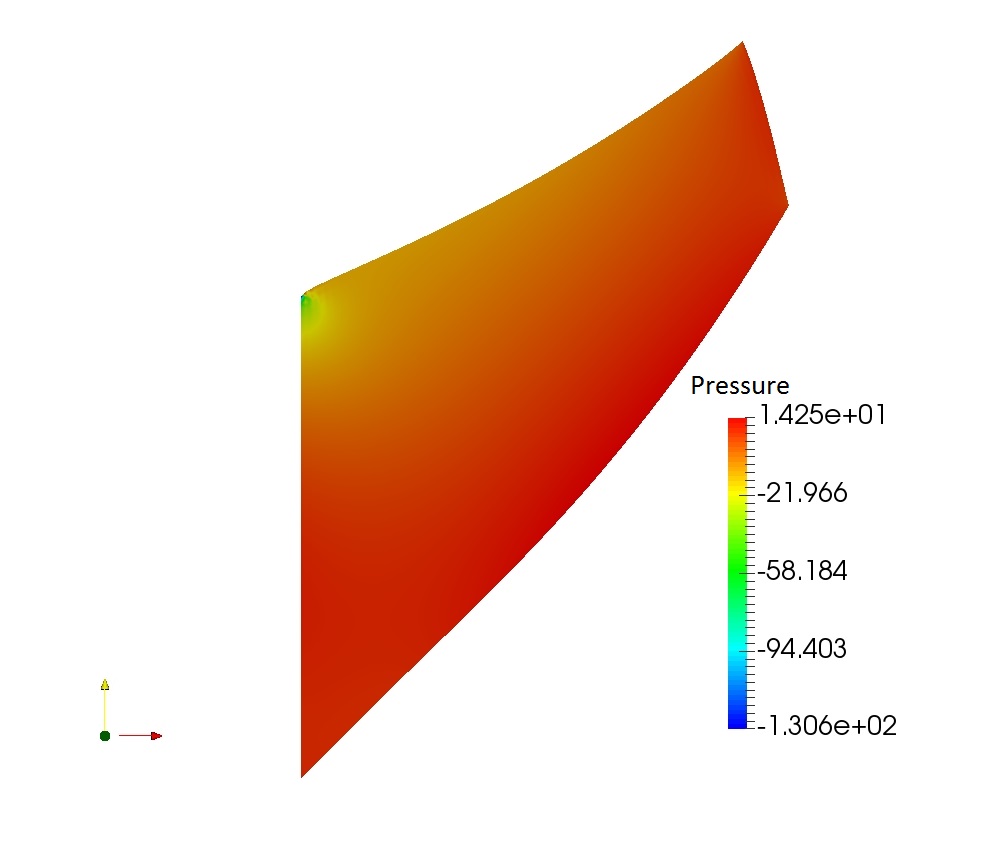}
    \caption{t=7s} \label{fig:2b}
  \end{subfigure}%

\caption{Pressure contours at different positions in time}
\label{fig:TransientCook}
\end{figure}

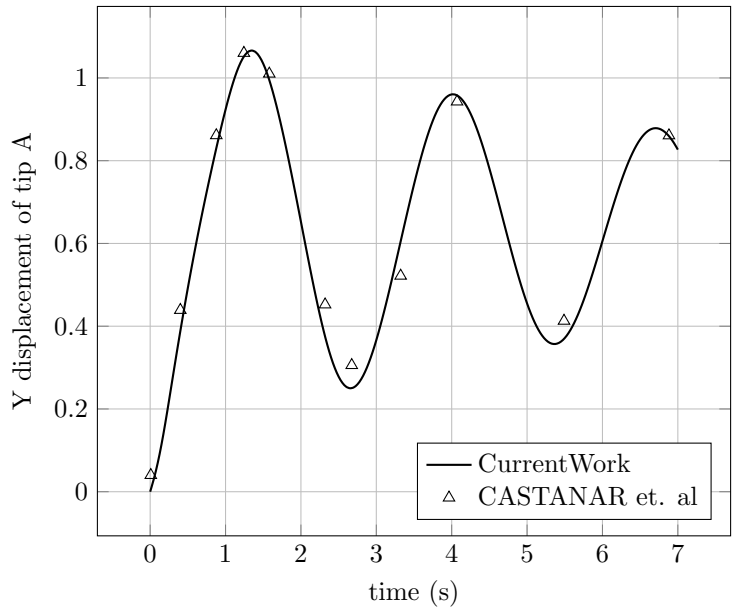
\begin{figure}
        \centering
        \begin{tikzpicture}[trim axis left, trim axis right]
                \begin{axis}[
                legend cell align=left,
                legend pos=south east,
                grid=major,
                xlabel={time (s)},
                ylabel={Y displacement of tip A}  
                ]                 
                
                \legend{CurrentWork, CASTANAR et. al}
                       
                \addplot[ 
                draw=black,
                thick,
                smooth
                ]
                table[
                x index=0,
                y index=1
                ]
                {data/YdispTransBDF1.dat};

                \addplot[ 
                only marks,
    		    style={solid, fill=gray},
    		    mark=triangle,
   			    mark size=2.5pt
                ]
                 table[
                x index=0,
                y index=1
                ]
                {data/YdispTransBenchmarkBDF1.dat};    
               
                \end{axis}
        \end{tikzpicture}
        \caption{Y displacement of tip A versus time}
		\label{fig:Variaiton}
\end{figure}

An extension of the aforementioned cook's membrane test to a transient regime is presented here.
The simulations are ran until t=7 s.
All the previous properties, along with the initial, and boundary conditions are preserved.
The geometry of the membrane was scaled with a factor of 0.1.
The solution obtained oscillates around the steady state solution computed earlier.
A graph showing the displacement of the tip A of the membrane with respect to time is show in figure [\ref{fig:Variaiton}].
Pressure contours at different positions in time can also be found in figure [\ref{fig:TransientCook}].

\subparagraph{Upsetting problem}

\begin{figure}
\centering
 \hspace*{\fill}   % maximize separation between the subfigures
\begin{subfigure}{0.45\textwidth}
   \includegraphics[width=1\linewidth]{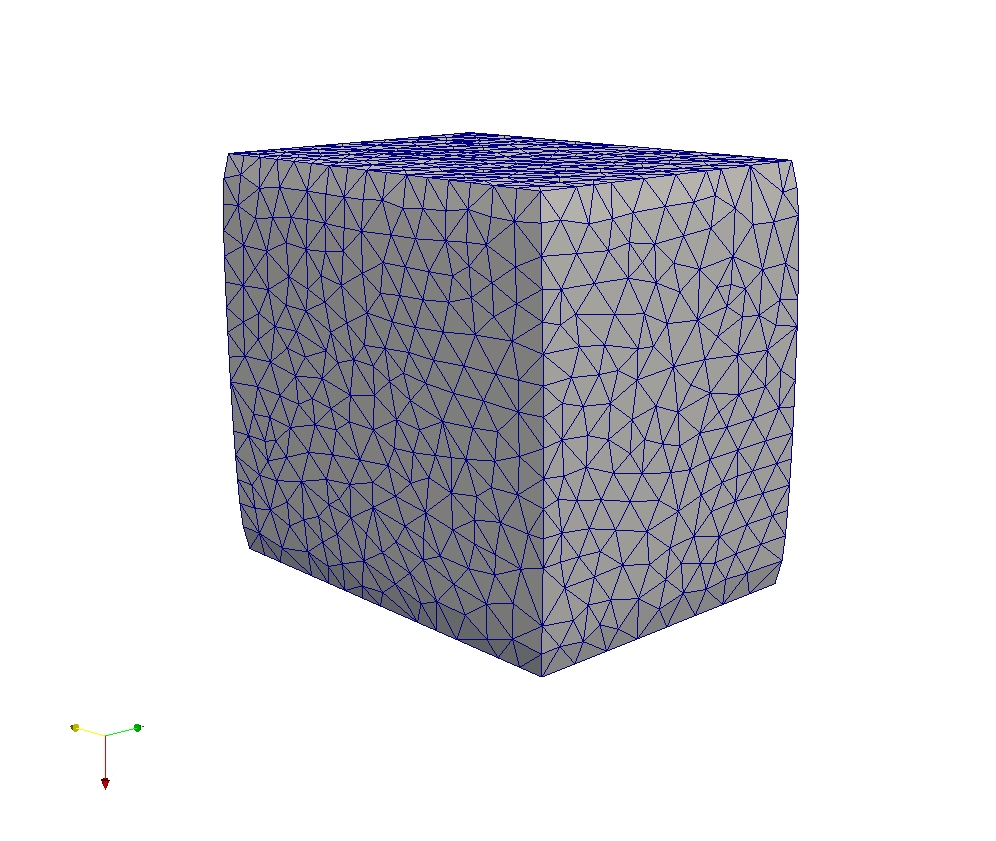}
    \caption{Mesh} \label{fig:4a}
  \end{subfigure}%
  \hspace*{\fill}   % maximize separation between the subfigures
  \begin{subfigure}{0.45\textwidth}
   \includegraphics[width=1\linewidth]{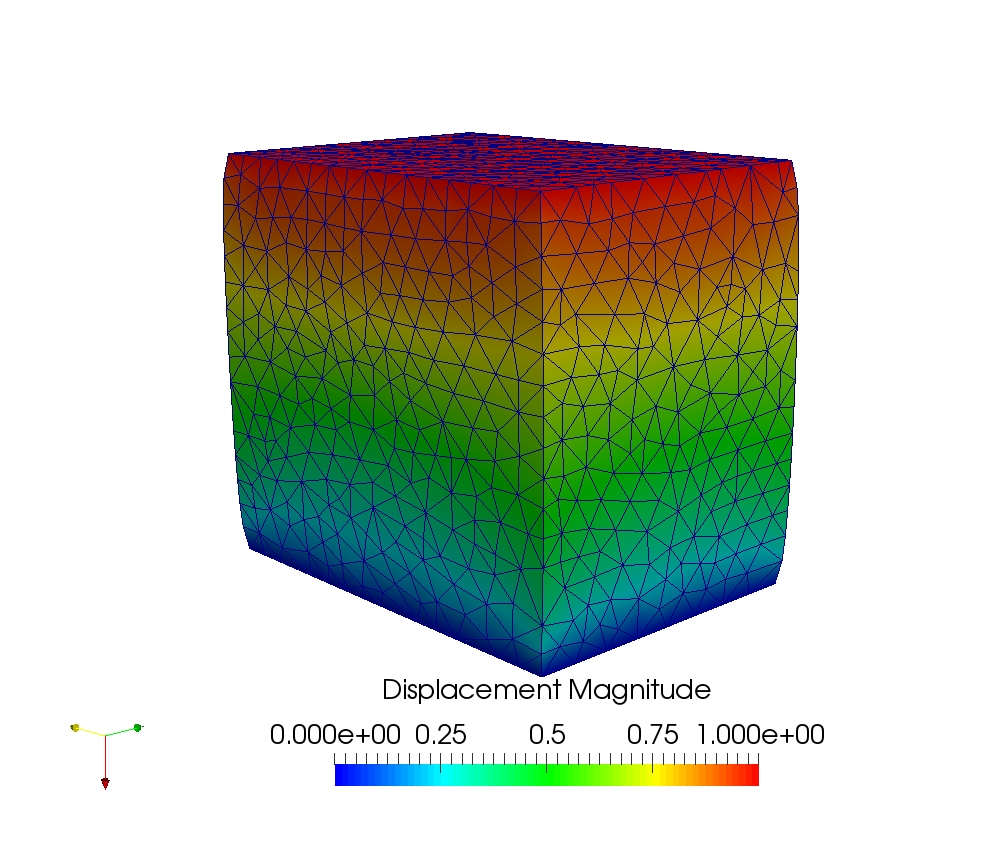}
    \caption{Displacement} \label{fig:4b}
  \end{subfigure}%
  \\
  \hspace*{\fill}   % maximize separation between the subfigures
  \begin{subfigure}{0.45\textwidth}
   \includegraphics[width=1\linewidth]{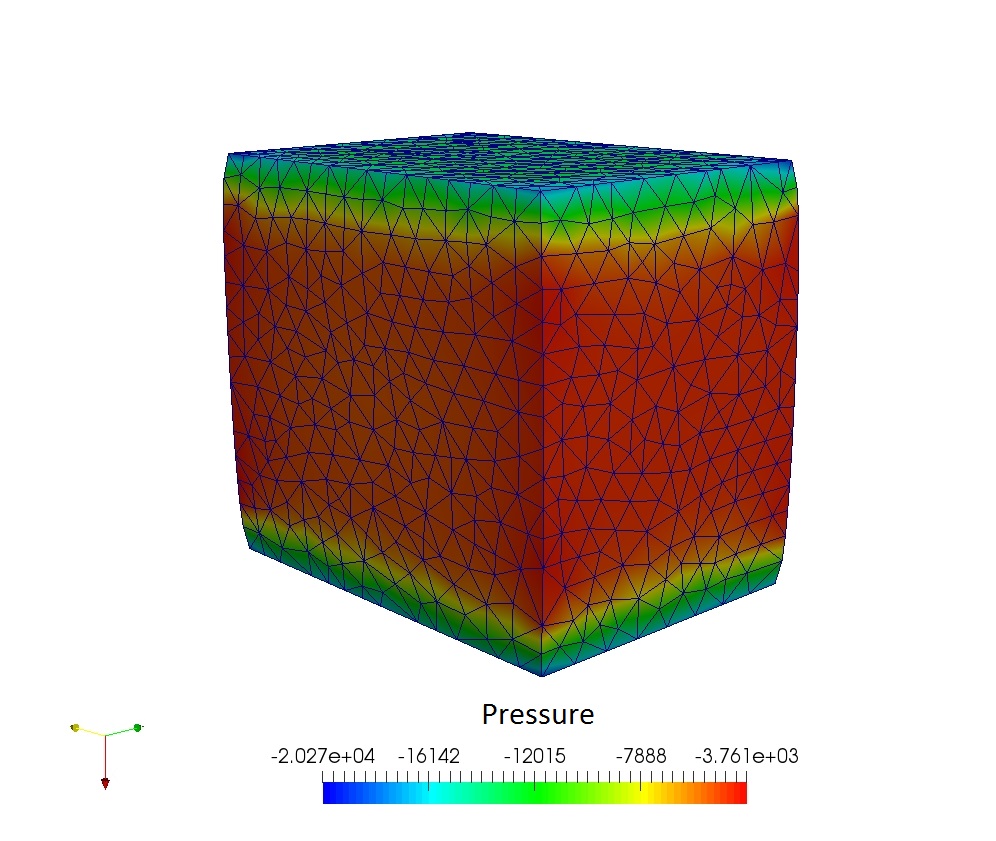}
    \caption{Pressure contours} \label{fig:4c}
  \end{subfigure}
  \hspace*{\fill}   % maximize separation between the subfigures
  \begin{subfigure}{0.45\textwidth}
   \includegraphics[width=1\linewidth]{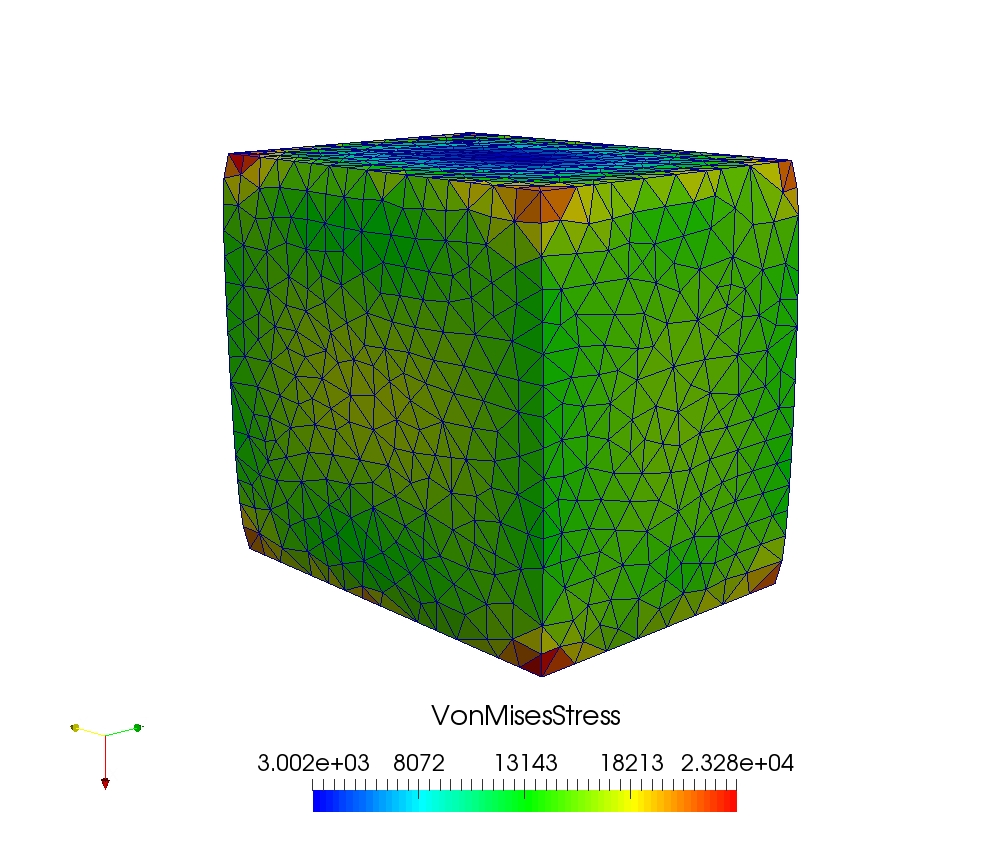}
    \caption{Von Mises Stress} \label{fig:4d}
  \end{subfigure}%

\caption{Upsetting problem results}
\label{fig:Upsetting}
\end{figure}

This example serves to show the capabilities of the formulation in the near incompressible case, even in the case of coarse tetrahedral meshes in 3D.
It consists of a 14x14x10 specimen that is deformed up to 7 $ \% $ of its height.
The specimen is clamped at the bottom, and clamped at the top to a rigid plate which gradually moves downwards while compressing the specimen.
Linear Elastic constitutive model is used.
The young's modulus and poisson ratio ar given as $E=2.0e+5 MPa$, and $\nu=0.4999$.
Figure \ref{fig:Upsetting} shows the behavior of our formulation.
For a relatively coarse unstructured 3D mesh, we can observe the displacement field, and the pressure contours.
The Von Mises stresses are also computed for the completeness of the study.
The results correlates with the literature and no pressure locking is observed.

\subsection{HyperElastic}

\subparagraph{A Computational Solid Mechanics test}

\begin{figure}
\centering
   \includegraphics[width=0.7\linewidth]{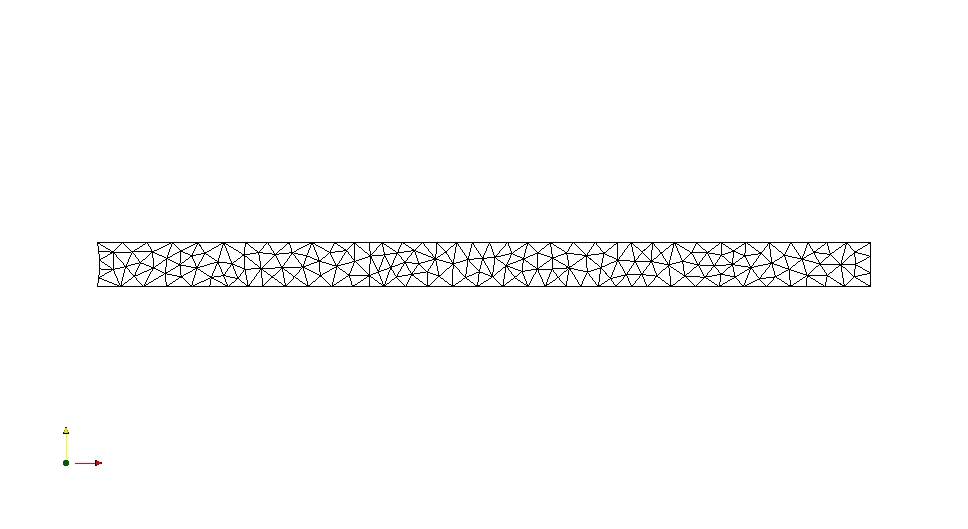}
\caption{Unstructured 2D mesh}
\label{fig:CSM2d}
\end{figure}

This test is a part of a well documented benchmark on Fluid--Structure Interaction, which deals with the solid part alone \cite{10.1007/3-540-34596-5_15}.
The structure is assumed to be elastic and compressible.
The constitutive law of the material is given by the St. Venant--Kirchhoff material.
The elastic beam is taken alone and subjected to a gravitational force $\gg=(0,g)$.
The beam has a length $l=0.35$, and a thickness $t=0.02$.
A typical unstructured mesh is shown in figre [\ref{fig:CSM2d}].
Three variations of the test are presented, two of which converges towards a steady-state solution and a non-diffusive transient case.
The different parameters are given in table [\ref{table:1}].
The results obtained are highly coherent with the literature.

\begin{table}[h!]
\centering
\begin{tabular}{ |p{3cm}||p{3cm}|p{3cm}|p{3cm}|  }
 \hline
 \multicolumn{4}{|c|}{CSM Cases} \\
 \hline
Solid Properties & Case 1 & Case 2 & Case 3\\
 \hline
$\rho$ & 1000 & 1000 & 1000\\
$nu$ & 0.4 & 0.4 & 0.4\\
$mu$ & 500000 & 2000000 & 500000\\
E & 1400000 & 5600000 &  1400000\\
g & 2 & 2 & 2\\
 \hline

\end{tabular} 

\caption{Different variations of the CSM test.}
\label{table:1}

\end{table}

\begin{figure}
\centering
  \includegraphics[width=0.8\linewidth]{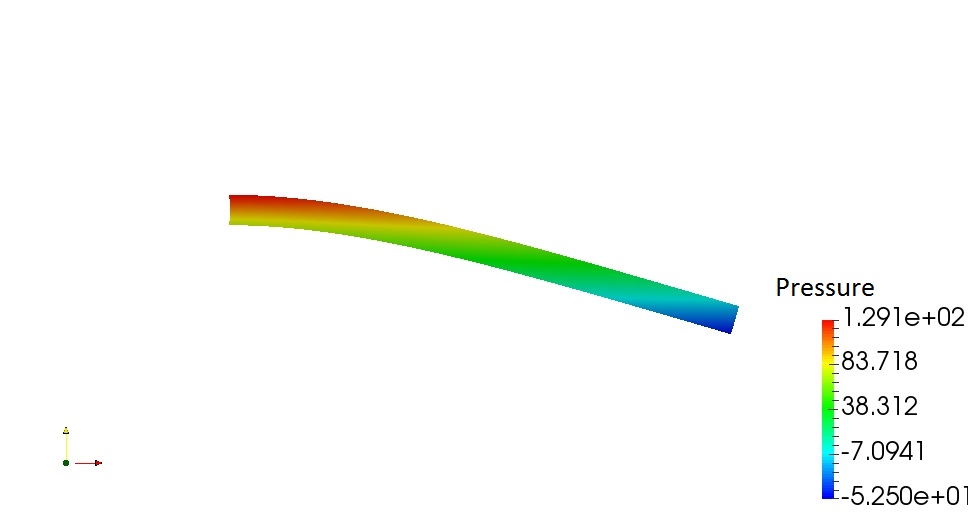}
\caption{Case 1} 
\label{fig:6a}
\end{figure}

\begin{figure}
\centering
  \includegraphics[width=0.8\linewidth]{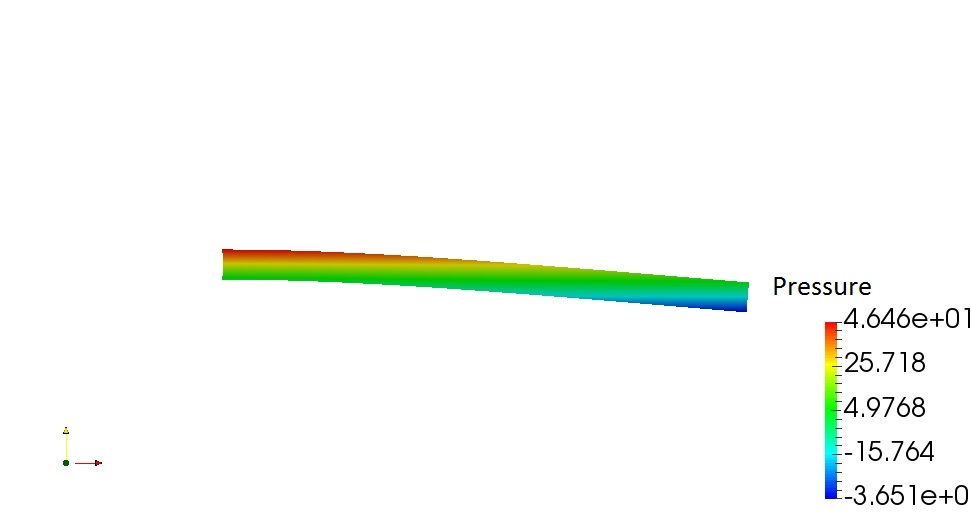}
\caption{Case 2} 
\label{fig:7b}
\end{figure}

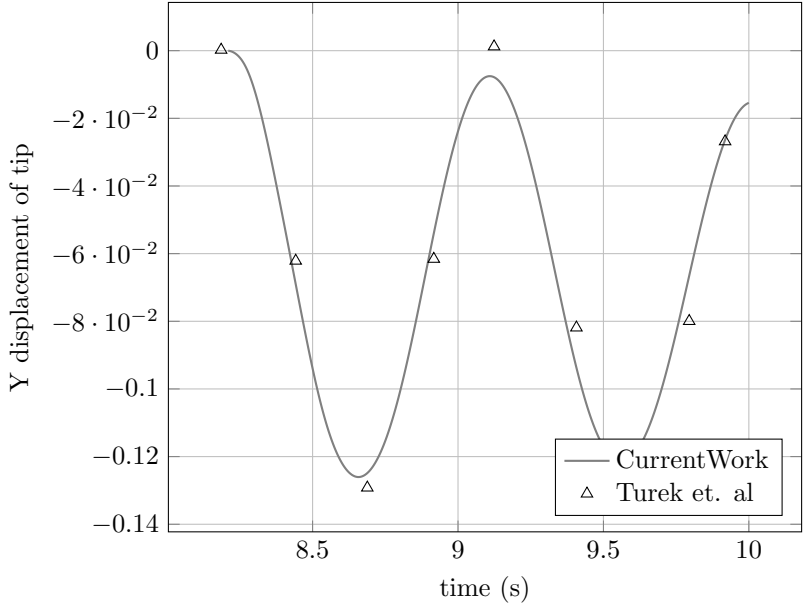
\begin{figure}
        \centering
        \begin{tikzpicture}[trim axis left, trim axis right]
                \begin{axis}[
                legend cell align=left,
                legend pos=south east,
                grid=major,
                xlabel={time (s)},
                ylabel={Y displacement of tip}
                ]                 
                
                \legend{CurrentWork, Turek et. al}  
                       
                \addplot[ 
                draw=gray,
                thick,
                smooth
                ]
                table[
                x index=0,
                y index=1
                ]
                {data/CSM.dat}; 
                
                \addplot[ 
                only marks,
    		    style={solid, fill=gray},
    		    mark=triangle,
   			    mark size=2.5pt
                ]
                 table[
                x index=0,
                y index=1
                ]
                {data/CSMBenchmark.dat};  
               
                \end{axis}
        \end{tikzpicture}
        \caption{Case 3: Y displacement versus time}
		\label{fig:CSM3}
\end{figure}

\subparagraph{Bending Beam test (3D)}

\begin{figure}
\centering
\begin{subfigure}{0.5\textwidth}
   \includegraphics[width=1\linewidth]{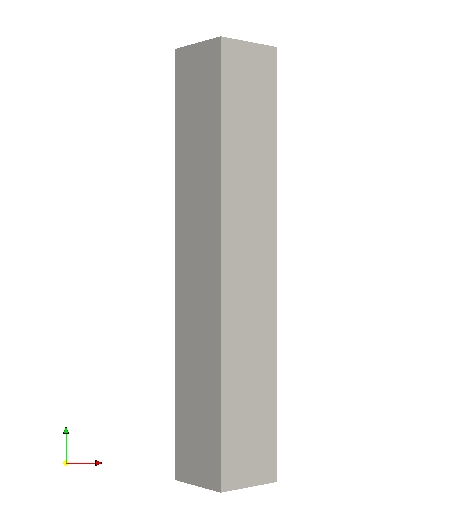}
    \caption{Geometrical Setup} \label{fig:5a}
  \end{subfigure}%
  \hspace*{\fill}   % maximize separation between the subfigures
  \begin{subfigure}{0.5\textwidth}
   \includegraphics[width=1\linewidth]{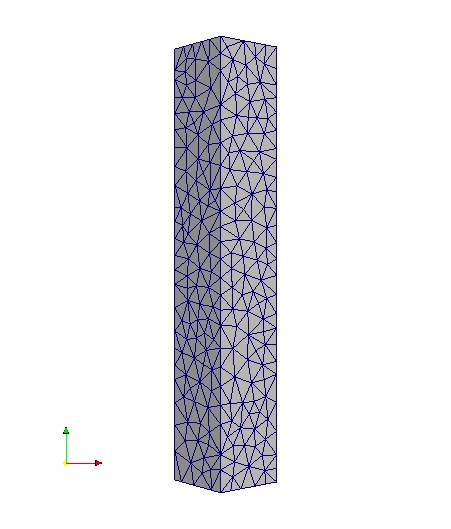}
    \caption{Unstructured mesh} \label{fig:5b}
  \end{subfigure}%
  \\
  \centering
  \begin{subfigure}{0.5\textwidth}
   \includegraphics[width=1\linewidth]{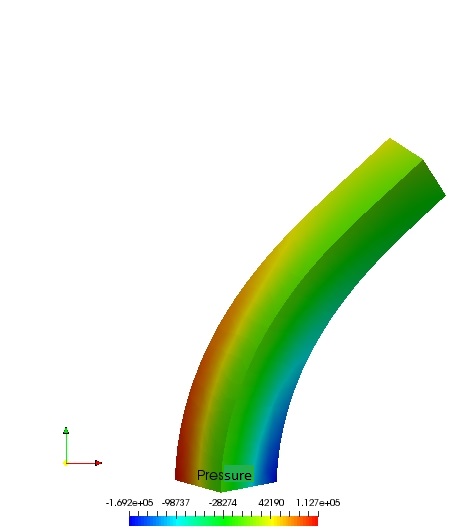}
    \caption{Pressure contours} \label{fig:5c}
  \end{subfigure}

\caption{3D Bending Beam Test}
\label{fig:BendingBeam}
\end{figure}

This test consists of a bending problem of a square cylinder in 3D.
The dimension of the square cylinder are given by 1x1x6 m.
The beam is also rotated with an angle of 5.2 degrees to avoid symmtery.
At $t=0$, the beam is stress free, and the displacement is equal to 0.
The geometry and a typical unstructured mesh is shown in figures [\ref{fig:5a}][\ref{fig:5b}] respectively.

An intial velocity is applied on the beam, given by:

\begin{equation}
\vv(\xx,0)=\vv(x,y,z,0)=(\frac{5y}{3},0,0)^T \: m/s \: y \in [0,6]m.
\end{equation}

The origin of our coordinates system is located at $(0.5,0.5,0)$.
The material is Neo-Hookean with the following properties:
$\rho_0=1.1e3$  $kg/m^3$, $E=1.7e7$  $Pa$, and $\nu=0.499$.

Zero displacement dirichlet boundary condition are imposed on the bottom of the cylinder, which is assumed to be clamped.
Zero traction boundary condition is applied on all other surfaces.

The simulation was ran until T= 2s.
This problem shows the capabilities of the framework in bending dominated problems of coarse unstructured mesh elements.
The unpolluted pressure field is shown in figure [\ref{fig:5c}]

\begin{figure}
\centering
 \hspace*{\fill}   % maximize separation between the subfigures
\begin{subfigure}{0.48\textwidth}
   \includegraphics[width=1\linewidth]{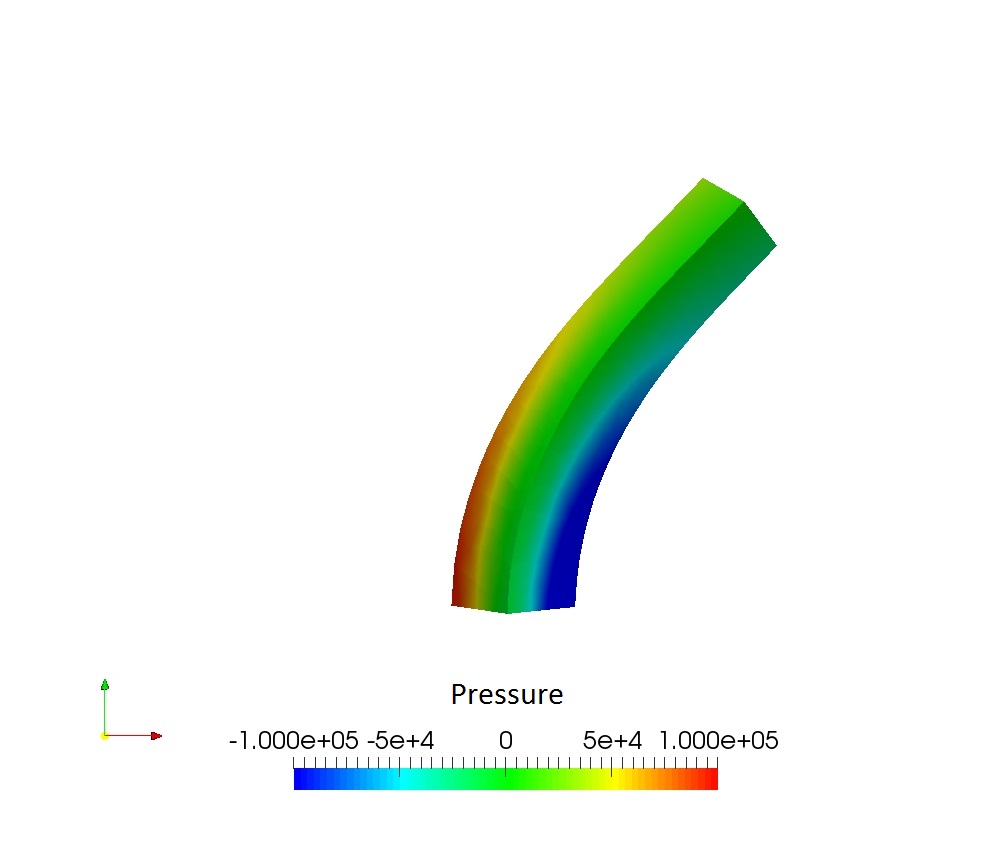}
    \caption{t=0.5s} \label{fig:7a}
  \end{subfigure}%
  \hspace*{\fill}   % maximize separation between the subfigures
  \begin{subfigure}{0.48\textwidth}
   \includegraphics[width=1\linewidth]{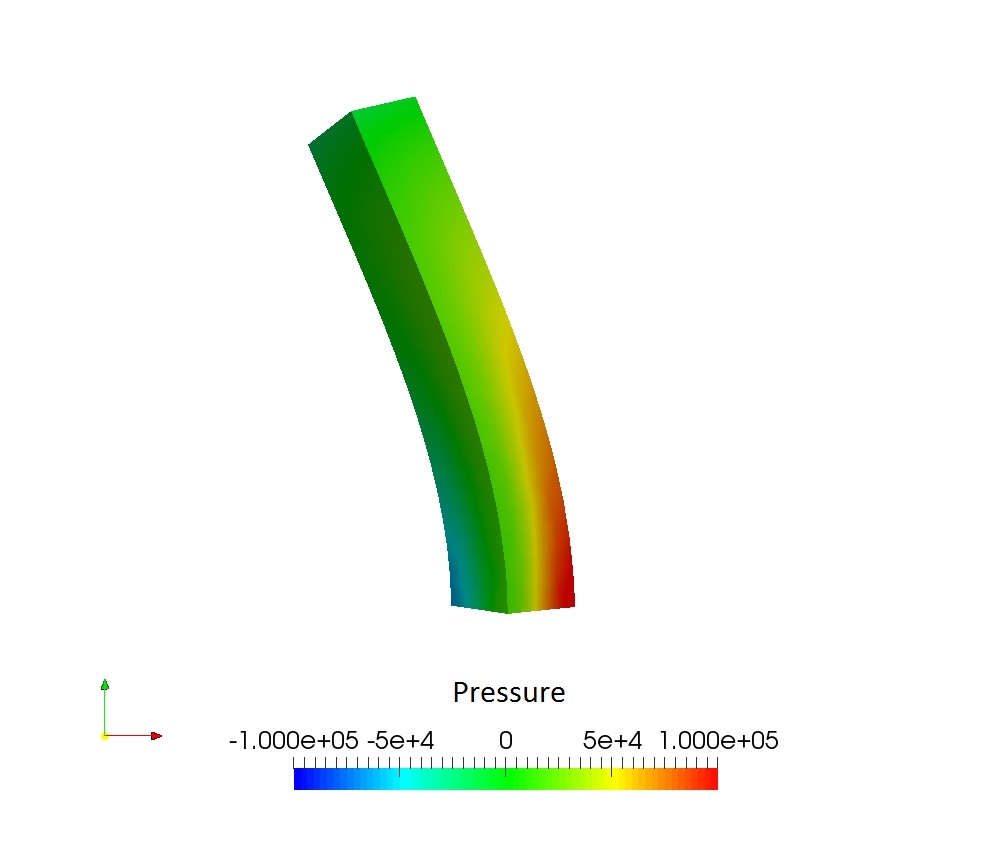}
    \caption{t=1s} \label{fig:7b}
  \end{subfigure}%
  \\
  \hspace*{\fill}   % maximize separation between the subfigures
  \begin{subfigure}{0.48\textwidth}
   \includegraphics[width=1\linewidth]{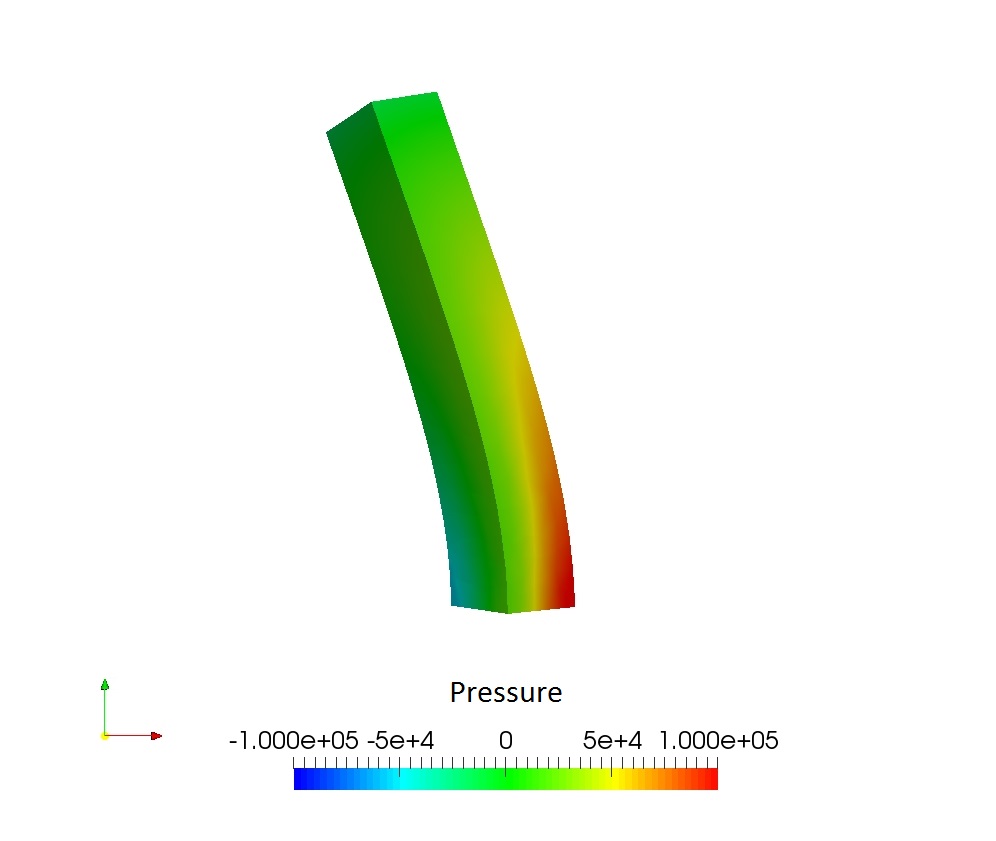}
    \caption{t=1.5s} \label{fig:7c}
  \end{subfigure}
  \hspace*{\fill}   % maximize separation between the subfigures
  \begin{subfigure}{0.48\textwidth}
   \includegraphics[width=1\linewidth]{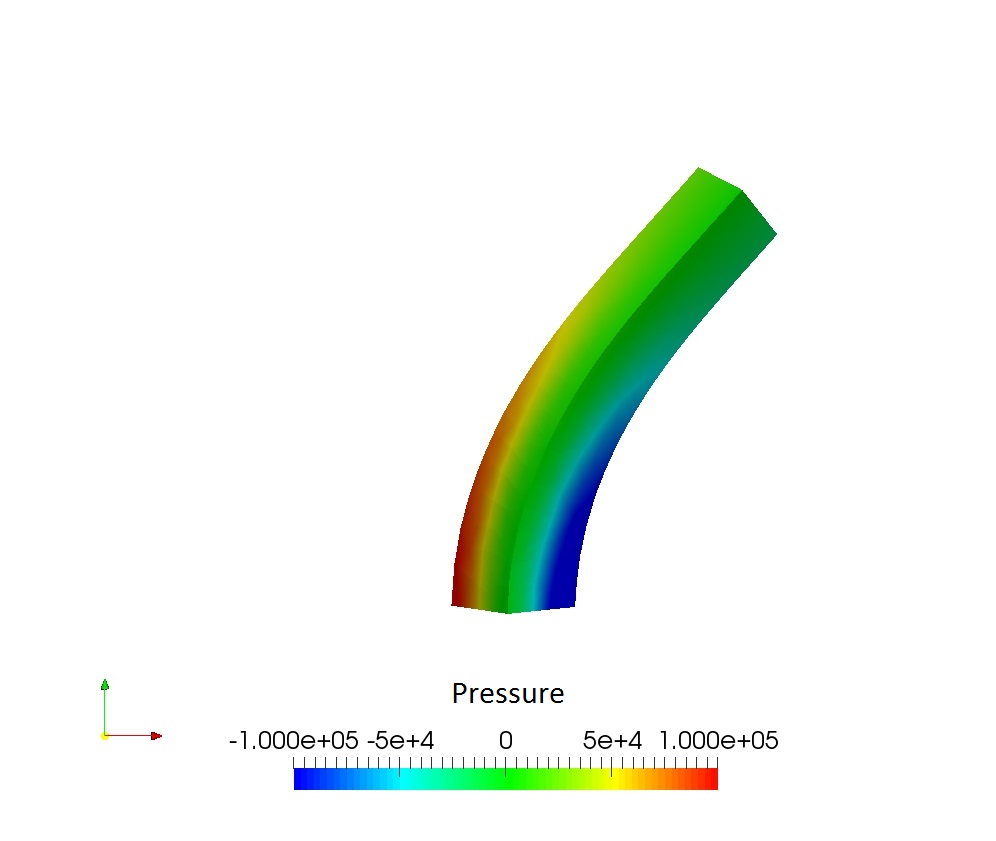}
    \caption{t=2s} \label{fig:7d}
  \end{subfigure}%

\caption{Upsetting problem results}
\label{fig:Transient3DBB}
\end{figure}

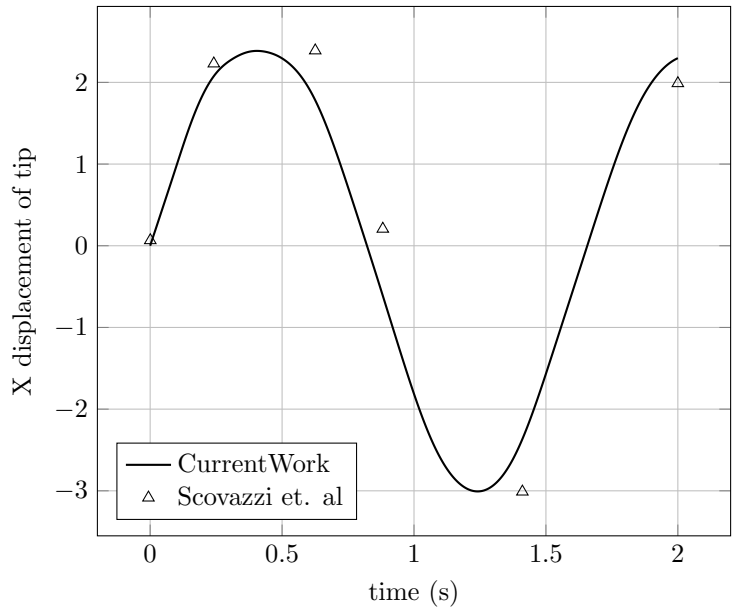
\begin{figure}
        \centering
        \begin{tikzpicture}[trim axis left, trim axis right]
                \begin{axis}[
                legend cell align=left,
                legend pos=south west,
                grid=major,
                xlabel={time (s)},
                ylabel={X displacement of tip}
                ]  
                
                \legend{CurrentWork, Scovazzi et. al}               
                       
                \addplot[ 
                draw=black,
                thick,
                smooth
                ]
                table[
                x index=0,
                y index=1
                ]
                {data/3DBB.dat};  
                
                \addplot[ 
                only marks,
    		    style={solid, fill=gray},
    		    mark=triangle,
   			    mark size=2.5pt
                ]
                 table[
                x index=0,
                y index=1
                ]
                {data/3DBBBenchmark.dat};  
               
                \end{axis}
        \end{tikzpicture}
        \caption{X displacement of tip versus time}
		\label{fig:XDisp}
\end{figure}

\subparagraph{Complex Geometry}

\begin{figure}
\centering
   \includegraphics[width=0.8\linewidth]{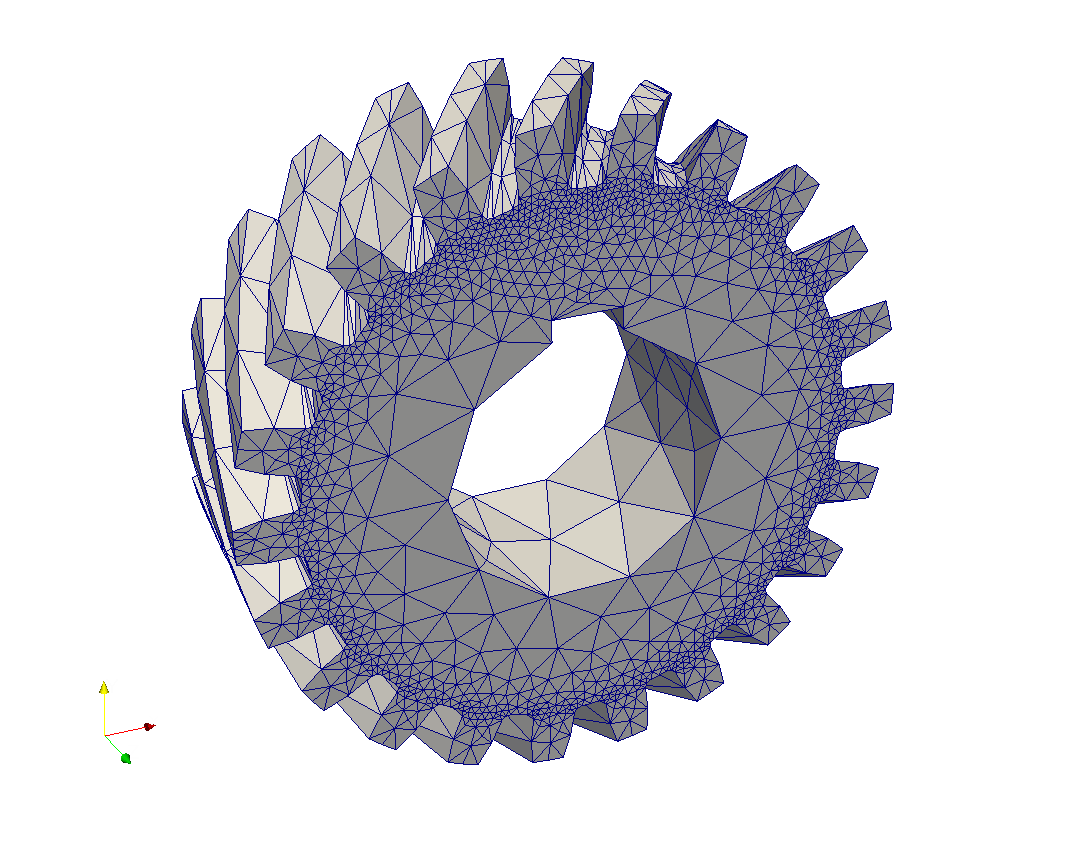}
    \caption{Geometrical Setup and computational mesh \cite{WinNT}} \label{fig:8}
\end{figure}

\begin{figure}
\centering
 \hspace*{\fill}   % maximize separation between the subfigures
\begin{subfigure}{0.45\textwidth}
   \includegraphics[width=1\linewidth]{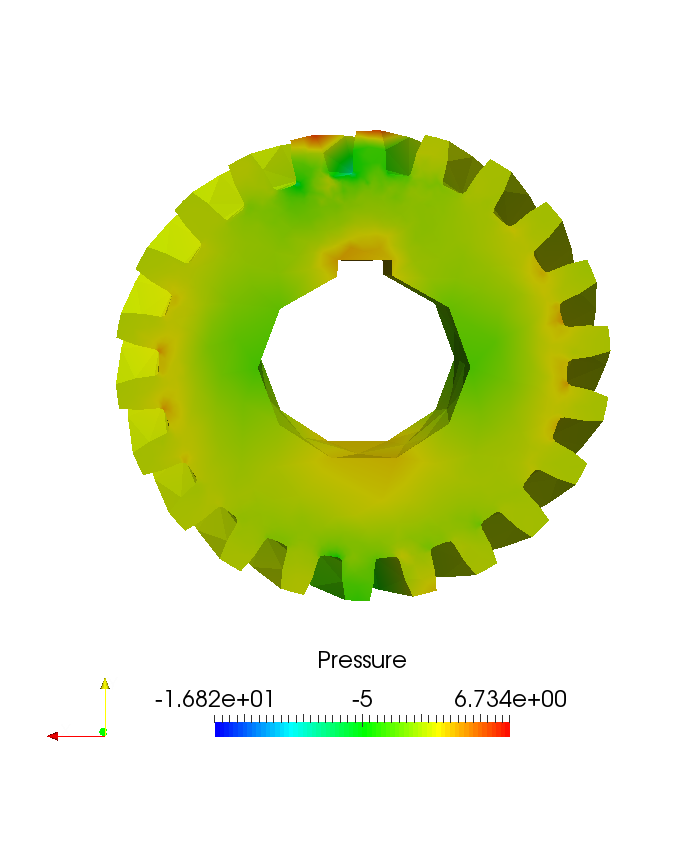}
    \caption{t=10s} \label{fig:9a}
  \end{subfigure}%
  \hspace*{\fill}   % maximize separation between the subfigures
  \begin{subfigure}{0.45\textwidth}
   \includegraphics[width=1\linewidth]{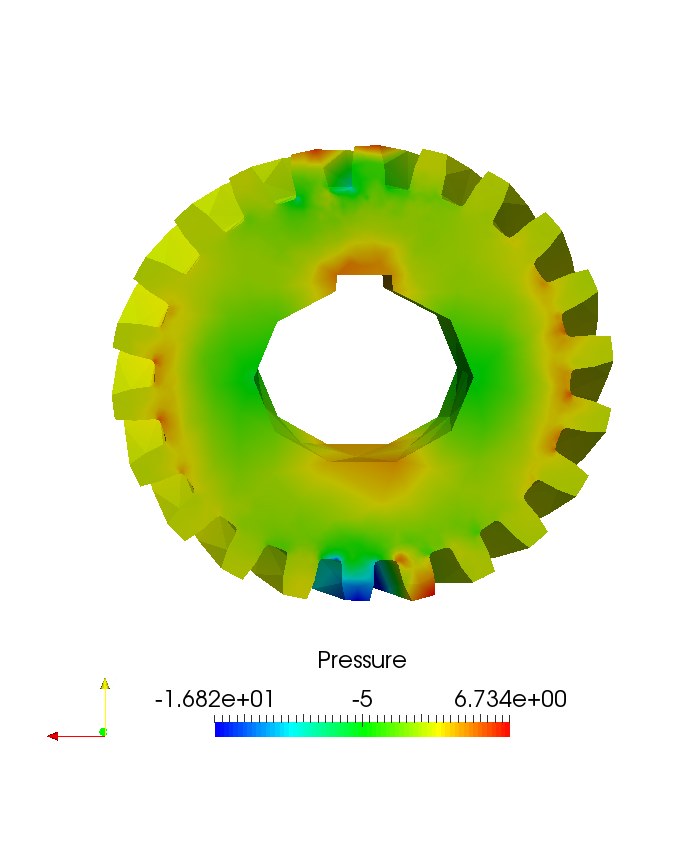}
    \caption{t=20s} \label{fig:9b}
  \end{subfigure}%
  \\
  \hspace*{\fill}   % maximize separation between the subfigures
  \begin{subfigure}{0.45\textwidth}
   \includegraphics[width=1\linewidth]{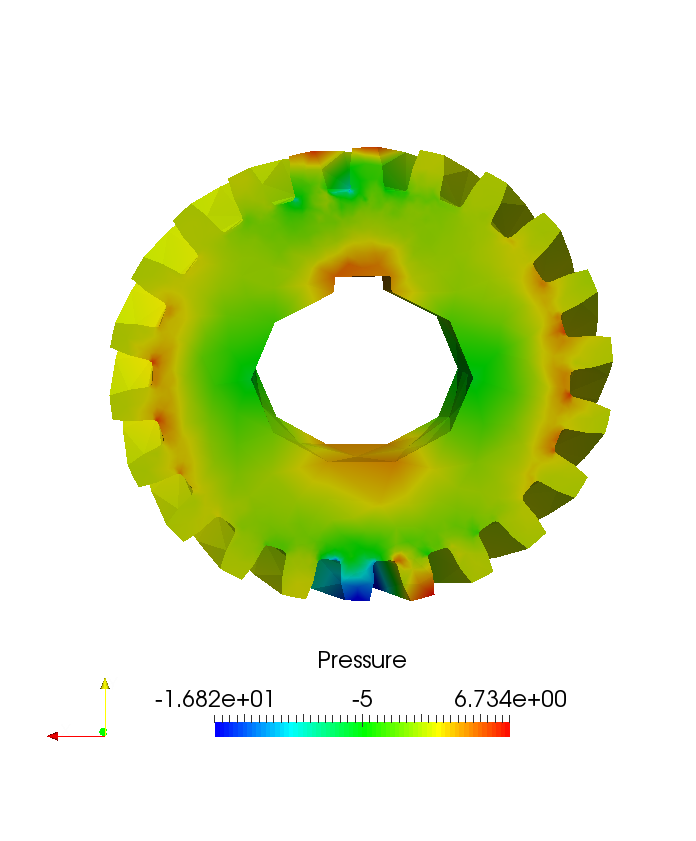}
    \caption{t=30s} \label{fig:9c}
  \end{subfigure}
  \hspace*{\fill}   % maximize separation between the subfigures
  \begin{subfigure}{0.45\textwidth}
   \includegraphics[width=1\linewidth]{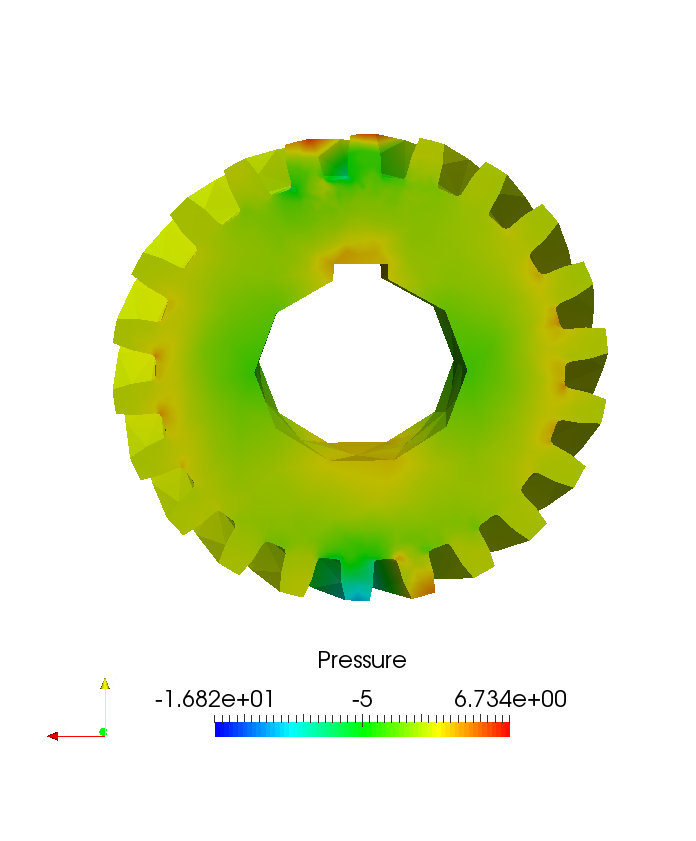}
    \caption{t=40s} \label{fig:9d}
  \end{subfigure}%
   \\
 %  \hspace*{\fill}   % maximize separation between the subfigures
  \begin{subfigure}{0.45\textwidth}
   \includegraphics[width=1\linewidth]{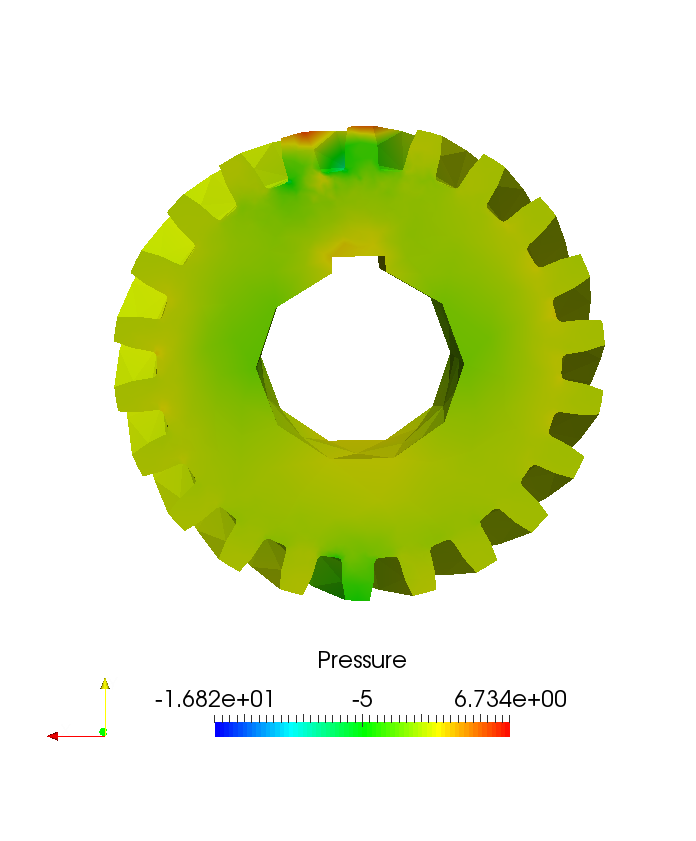}
    \caption{t=50s} \label{fig:9e}
  \end{subfigure}

\caption{Pressure Contours at different time steps}
\label{fig:TransientGear}
\end{figure}

As a final numerical test, the ability of the framework to handle complex geometries is evaluated.
A helical gear is shown in figure [\ref{fig:8}], along with its computational mesh, consisting of 9865 tetrahedral 3D elements.
The material is assumed to have a non-linear transient elastic behavior, and is considered fully incompressible with $\rho=1$, $E=250$, and $\nu=0.5$.
The Neo-Hookean model is used.
A vertical downward forcing term of magnitude 5 is imposed on the top plane of the geometry.
Homogeneous Dirichlet boundary conditions are imposed on the bottom plane of the geometry.
Zero traction neumann boundary conditions are applied on the rest of the boundaries.
Pressure contours at times equal to 10,20,30,40, and 50 s are shown in figure [\ref{fig:TransientGear}].
High pressure gradients can be observed around the hole of the geometry.
The solution converged properly, and no spurious pressure oscillation was observed.

\section{Perspectives and conclusion}
\indent

In this work, we have presented a framework based on unstructured tetrahedral meshes that can handle complex geometries, which models the nonlinear behavior of solid elastodynamics.
By combining the proposed new mixed formulation in the updated Lagrangian framework, and the R-method for moving meshes, the framework was able to handle nearly, and fully incompressible material in bending dominated problems.
This was achieved through the deviatoric/volumetric split of the stress tensor.
A piece wise linear mixed formulation in displacement and pressure was obtained in the updated Lagrangian formulation.
We achieved piece wise linear interpolation for both displacement and pressure through the Variational Multi-Scale approach, based on the orthogonal decomposition of the function spaces (extended from fluid mechanics).
The stabilization proved effective in both steady-state and transient regime.
We are in the process of applying this newly developed solver for Fluid--Structure Interaction (FSI) applications.

%===============================================================================

\clearpage
\biboptions{sort&compress}
\bibliographystyle{elsarticle-num}
%\thebibliography{biblio}
\bibliography{SolidModeling}
\end{document}